\documentclass[a4paper, 11pt, reqno]{amsart}

\usepackage{amsmath, latexsym, amsfonts, amssymb, amsthm, amscd,
  stmaryrd}

\usepackage[T1]{fontenc}

\usepackage[english]{babel}

\usepackage{graphicx}

\numberwithin{equation}{section}

\newtheorem{theorem}{Theorem} 
\newtheorem{lemma}{Lemma}

\newtheorem{corollary}{Corollary}

\theoremstyle{remark}
\newtheorem*{remark}{Remark}

\theoremstyle{definition}
\newtheorem{definition}{Definition}

\newcounter{assumption}

\renewcommand{\leq}{\leqslant}
\renewcommand{\geq}{\geqslant}

\newcounter{hyp}

\newcommand{\etal}{\textit{et al. }}

\newcommand{\mc}{\mathcal}
\newcommand{\mb}{\mathbf}
\newcommand{\mbb}{\mathbb}
\newcommand{\mrm}{\mathrm}
\newcommand{\mf}{\mathfrak}
\newcommand{\mbt}[1]{\wt{\mb #1}}

\newcommand{\tup}[1]{\textup{#1}}

\newcommand{\ellom}{\ell_{\omega}}

\newcommand{\wh}{\widehat}

\newcommand{\wt}{\widetilde}

\newcommand{\kpa}{\kappa}

\newcommand{\von}{\varepsilon}

\newcommand{\tta}{\theta}

\newcommand{\Lba}{\Lambda}

\newcommand{\lba}{\lambda}

\newcommand{\raro}{\rightarrow}

\newcommand{\laro}{\leftarrow}

\newcommand{\setR}{\mathbb{R}}

\newcommand{\setN}{\mathbb{N}}

\newcommand{\setRp}{\setR^{+}}

\newcommand{\RV}{\mrm{RV}}

\newcommand{\agrid}{\mc H_a^{\textup{arith}}}

\newcommand{\ggrid}{\mc H_a^{\textup{geom}}}

\newcommand{\eqdef}{\triangleq}

\newcommand{\ind}[1]{\mathbf{1}_{#1}}

\DeclareMathOperator*{\argmax}{argmax}

\newcommand{\prodsca}[2]{\langle #1 \, ,\, #2 \rangle}

\newcommand{\prodscah}[2]{\prodsca{#1}{#2}_{h}}
\newcommand{\prodscahp}[2]{\prodsca{#1}{#2}_{h'}}

\newcommand{\norm}[1]{\| #1 \|}

\newcommand{\norminfty}[1]{\| #1 \|_{\infty}}

\newcommand{\normh}[1]{\norm{#1}_{h}}
\newcommand{\normhp}[1]{\norm{#1}_{h'}}

\newcommand{\trans}{{}^{t}}

\newcommand{\ppint}[1]{\lfloor #1 \rfloor} 

\newcommand{\Prob}{\mbb{P}}

\newcommand{\Pm}{\Prob_{\mu}^{n}}
\newcommand{\Pfm}{\Prob_{f,\mu}^{n}}
\newcommand{\Efm}{\mbb E_{\,f, \mu}^{\,n}}

\newcommand{\Em}{\mbb E_{\mu}^{n}}

\newcommand{\Var}{\mbb V \text{ar}}

\newlength{\figurelength} 
\title[Pointwise adaptive curve estimation with degenerate design]{On
  pointwise adaptive curve estimation with a degenerate random design}

\author{St{\'e}phane Ga{\"\i}ffas}

\address{Laboratoire de Probabilit{\'e}s et Mod{\`e}les
  Al{\'e}atoires, U.M.R. CNRS 7599 and Universit{\'e} Paris 7, 175 rue
  du Chevaleret, 75013 Paris}

\email{gaiffas@math.jussieu.fr}

\keywords{adaptive estimation, degenerate design, nonparametric
  regression, random design. }

\subjclass[2000]{62G05, 62G08}

\date{\today}

\dedicatory{Laboratoire de Probabilit{\'e}s et Mod{\`e}les
  Al{\'e}atoires \\
  Universit{\'e} Paris 7, 175 rue du Chevaleret, 75013 Paris \\
  email: \rm{ \texttt{gaiffas@math.jussieu.fr} } }

\begin{document}

\selectlanguage{english}

\begin{abstract}
  We consider the nonparametric regression with a random design model,
  and we are interested in the adaptive estimation of the regression
  at a point $x_0$ where the design is degenerate. When the design
  density is $\beta$-regularly varying at $x_0$ and $f$ has a
  smoothness $s$ in the H\"older sense, we know from Ga\"iffas
  (2004)\nocite{gaiffas04a} that the minimax rate is equal to
  $n^{-s/(1+2s+\beta)} \ell(1/n)$ where $\ell$ is slowly varying. In
  this paper we provide an estimator which is adaptive both on the
  design and the regression function smoothness and we show that it
  converges with the rate $(\log n/n)^{s/(1+2s+\beta)} \ell(\log
  n/n)$. The procedure consists of a local polynomial estimator with a
  Lepski type data-driven bandwidth selector similar to the one in
  Goldenshluger and Nemirovski
  (1997)\nocite{goldenshluger_nemirovski97} or Spokoiny
  (1998)\nocite{spok98}. Moreover, we prove that the payment of a
  $\log$ in this adaptive rate compared to the minimax rate is
  unavoidable.
\end{abstract}

\maketitle

\section{Introduction}
\label{sec:introduction}

\subsection{The model}
\label{sec:the_model}

We observe $n$ pairs of random variables $(X_i, Y_i) \in \setR \times
\setR$ independent and identically distributed satisfying
\begin{equation}
  \label{eq:regression_model}
  Y_i = f(X_i) + \xi_i,
\end{equation}
where $f : [0,1] \rightarrow \setR$ is the unknown signal to be
recovered, the variables $(\xi_i)$ are centered Gaussian with variance
$\sigma^2$ and independent of the design $X_1, \ldots, X_n$. The
variables $X_i$ are distributed with respect to a density $\mu$. We
want to recover $f$ at a fixed point $x_0$.

The classical way to consider the nonparametric regression model is to
take $X_i = i / n$. In this model with an equispaced design the
observations are \emph{homogeneously} distributed over the unit
interval. If we take the $X_i$ random we can modelize cases with
\emph{inhomogeneous} observations as the design distribution is "far"
from the uniform law. We allow here the density $\mu$ to be
\emph{degenerate} (vanishing or exploding) and we are more precisely
interested in the adaptive estimation of $f$ at a point where the
design is degenerate, namely a point with very inhomogeneous data.

\subsection{Motivations}
\label{sec:motivations}

The adaptive estimation of the regression function is a well-developed
problem. Several adaptive procedures can be applied for the estimation
of a function with unknown smoothness: nonlinear wavelet estimation
(thresholding), model selection, kernel estimation with a variable
bandwidth (the Lepski method), and so on.

Recent results dealing with the adaptive estimation of the regression
function when the design is not equispaced or random include
Antoniadis \etal (1997)\nocite{antoniadis_et_al97}, Brown and Cai
(1998)\nocite{cai_brown98}, Wong and Zheng
(2002)\nocite{wong_zheng02}, Maxim (2003)\nocite{voichitaphd},
Delouille \etal (2004)\nocite{delouille_et_al04}, Kerkyacharian and
Picard (2004)\nocite{kerk_picard_warped_04}, among others. A natural
question arises: what happens if we want to estimate adaptatively the
regression function at a point where the design is degenerate? In
Ga\"iffas (2004)\nocite{gaiffas04a} we proved when $\mu$ varies
regularly at $x_0$ that the minimax convergence rate $\psi_n$ over a
H{\"o}lder type regularity class with smoothness $s>0$ (around $x_0$)
satisfies
\begin{equation*}
  \psi_n \asymp n^{-s/(1+2s+\beta)} \ell(1/n) \text{ as } n \raro
  +\infty,
\end{equation*}
where $\beta$ is the regular variation index of $\mu$ at $x_0$ (see
definition \ref{def:rv_def}) and $\ell$ is slowly varying (the
notation $a_n \asymp b_n$ means $0 < \liminf a_n / b_n \leq \limsup
a_n / b_n < +\infty$). For the proof of the upper bound, a (non
adaptive) linear procedure was used.

The next logical step is then to find a procedure able to recover $f$
with as less prior knowledge as possible on its smoothness and on the
design density. On pointwise adaptive curve estimation (in the
regression or the white noise model) see Lepski
(1990)\nocite{lepski90}, Lepski and Spokoiny
(1997)\nocite{lepski_spok97}, Spokoiny (1998)\nocite{spok98} and Brown
and Cai (1998)\nocite{cai_brown98} for wavelet methods.

\subsection{Organisation of the paper}

We introduce the estimator in section \ref{sec:the_procedure}. In
section \ref{sec:upper_bounds} we give upper bounds for this procedure
conditionally on the design, see theorem
\ref{thm:NA_adapt_upper_bound} and in the regular variation framework,
see theorem \ref{thm:asympt_adaptive_upper_bound}. In section
\ref{sec:optimality} we prove that the obtained convergence rate is
optimal, see theorem \ref{thm:lower_bound} and its corollary.  We
present numerical illustrations in section \ref{sec:simulations} for
several datasets and we discuss in detail some points in section
\ref{sec:discussion}. Section \ref{sec:proofs} is devoted to the
proofs and we recall some well-known facts on regularly varying
functions in appendix.

\section{The procedure}
\label{sec:the_procedure}

\subsection{Local polynomial estimation}
\label{sec:local_polynomial_estimation}

Let $\kpa \in \setN$ and $h > 0$ (the \emph{bandwidth}). We define
\begin{equation*}
  N_{n,h} \eqdef \#\{ X_i \text{ such that } X_i \in [x_0-h,x_0+h] \},
\end{equation*}
and we introduce the pseudo-scalar product
\begin{equation*}
  \prodsca{f}{g}_h \eqdef \frac{1}{N_{n,h}} \sum_{|X_i-x_0| \leq h}
  f(X_i) g(X_i),
\end{equation*}
and $\norm{\cdot}_h$ the corresponding pseudo-norm. Let $\phi_{j}(x) =
(x-x_0)^j$ for $j = 0, \ldots, \kpa$. We introduce the matrix $\mb
X_h$ and the vector $\mb Y_h$ with entries for $0 \leq j,l \leq \kpa$:
\begin{equation}
  \label{eq:matrix_X_vector_Y_def}
  (\mb X_h)_{j,l} = \prodscah{\phi_{j}}{\phi_{l}} \quad \text{ and }
  \quad (\mb Y_h)_j = \prodscah{Y}{\phi_{j}}. 
\end{equation}
\begin{definition}
  \label{def:loc_pol_est_final}
  Let
  \begin{equation*}
    \wh f_{h,\kpa} =
    \begin{cases}
      \, \wh \tta_{h, 0} \phi_{0} + \wh \tta_{h, 1} \phi_{1} + \cdots
      + \wh \tta_{h, \kpa} \phi_{\kpa} &\text{ when } N_{n,h} > 0, \\
      \, 0 &\text{ when } N_{n,h} = 0,
    \end{cases}
  \end{equation*}
  where $\wh \tta_h$ is the solution of the linear system
  \begin{equation}
    \label{eq:lin_sys_final}
    \mbt X_h \tta = \mb Y_h,
  \end{equation}
  where
  \begin{equation*}
    \mbt X_h \eqdef \mb X_h + N_{n, h}^{-1/2} \mb I_{\kpa+1}
    \ind{\lba(\mb X_h) \leq N_{n,h}^{- 1/2}},
  \end{equation*}
  with $\lba(M)$ standing for the smallest eigenvalue of a matrix $M$
  and $\mb I_{\kpa+1}$ the identity matrix in $\setR^{\kpa+1}$. 
\end{definition}

This procedure is slightly different from the classical version of the
local polynomial estimator. We note that the correction term in $\mbt
X_h$ entails $\lba(\mbt X_h) \geq N_{n,h}^{-1/2}$. On local polynomial
estimation, see Stone (1980)\nocite{stone80}, Fan and Gijbels (1995,
1996)\nocite{fan_gijbels95}\nocite{fan_gijbels96}, Spokoiny
(1998)\nocite{spok98} and Tsybakov (2003)\nocite{tsybakov03} among
many others.

\subsection{Adaptive bandwidth selection}
\label{sec:the_bandwidth_selector}

The procedure selects the bandwidth $h$ in a set $\mc H$ called the
\emph{grid}, which is a tuning parameter of the adaptive procedure. We
can choose either an \emph{arithmetical} or a \emph{geometrical} grid
\begin{equation*}
  \mc H = 
  \begin{cases}
    \displaystyle \mc H_{a}^{\text{arith}} = \bigcup_{i=1}^{[(n-2)/a]}
    \{ h_{2 + [i a]} \} &\text{ for } a \geq 1, \text{ or } \\
    \displaystyle \mc H_a^{\text{geom}} = \bigcup_{i=1}^{[\log_a{n}]}
    \{ h_{[a^i]} \} &\text{ for } a > 1,
  \end{cases}
\end{equation*}
where $h_i \eqdef |X_{(i)} - x_0|$ and where $|X_{(i)} - x_0| \leq
|X_{(i+1)} - x_0|$ for any $i = 1, \ldots, n - 1$. Note that $[x]$
stands for the integer part of $x$. We define
\begin{equation*}
  \mc H_h \eqdef \{ h^{'} \in \mc H \text{ such that } h^{'} \leq h
  \}.
\end{equation*}
The bandwidth is selected as follows:

\begin{equation}
  \label{eq:H_n_hat_def}
  \wh H_n \eqdef \max\Bigl\{ h \in \mc H \big| \forall h^{'} \in \mc
  H_h \, \forall 0 \leq j \leq \kpa,\, |\prodsca{\wh f_{h,\kpa} - \wh
    f_{h',\kpa}}{{\phi_j}}_{h'}| \leq \sigma \norm{\phi_j}_{h'} T_{n,h',h}
  \Bigr\},
\end{equation}
where $\wh f_{h,\kpa}$ is given by definition
\ref{def:loc_pol_est_final} and where the threshold $T_{n,h',h}$ is
equal to
\begin{equation}
  \label{eq:T_n_def}
  T_{n,h',h} \eqdef
  \begin{cases}
    \displaystyle C_{\kpa} \sqrt{ C_{p} N_{n,h'}^{-1} \log N_{n,h}} +
    \sqrt{(N_{n,h} - a)^{-1} \log n} &\text{ if } \mc H =
    \agrid, \\
    C_{\kpa} \sqrt{C_{p} N_{n,h'}^{-1} \log N_{n,h}} + \sqrt{(1+a)
      N_{n,h}^{-1} \log n} &\text{ if } \mc H = \ggrid,
  \end{cases}
\end{equation}
with $C_{\kpa} \eqdef 1 + \sqrt{\kpa+1}$, $C_{p} = 8(1 + 2 p)$ where
$p$ fits with the loss function in \eqref{eq:def_maximal_risk} and $a$
is the grid parameter. The estimator is then
\begin{equation}
  \label{eq:final_estimator_definition}
  \wh f_n(x_0) \eqdef \wh f_{\wh H_n, \kpa}(x_0). 
\end{equation}
The selection rule \eqref{eq:H_n_hat_def} is similar to the method by
Lepski, see Lepski (1990)\nocite{lepski90}, Lepski \etal
(1997)\nocite{lepski_mammen_spok_97} and Lepski and Spokoiny
(1997)\nocite{lepski_spok97} and is additionally to the original
Lepski method sensitive to the design. This procedure is close to the
one in Spokoiny (1998)\nocite{spok98}. See section
\ref{sec:comparison_with_existing_results} for more details on
existing procedures in the literature.

\section{Upper bounds}
\label{sec:upper_bounds}

We measure a procedure $\wt f_n$ performance over a class $\Sigma$ (to
be specified in the following) with the maximal risk
\begin{equation}
  \label{eq:def_maximal_risk}
  \bigl(\sup_{\displaystyle f \in \Sigma}\Efm \{ | \wt f_n(x_0) -
  f(x_0)|^p \}\bigr)^{1/p},
\end{equation}
where $x_0$ is the estimation point and $p \geq 1$. The expectation
$\Efm$ in \eqref{eq:def_maximal_risk} is taken with respect to the
joint law $\Pfm$ of the observations \eqref{eq:regression_model}.

\subsection{Regular variation}
\label{sec:regular_variation_intro}

The regular variation definition and main properties are due to
Karamata (1930). On this topic we refer to Senata
(1976)\nocite{senata76}, Geluk and de Haan
(1987)\nocite{geluk_de_haan87}, Resnick (1987)\nocite{resnick87} and
Bingham \etal (1989)\nocite{bgt89}.

\begin{definition}[Regular variation]
  \label{def:rv_def}
  A continuous function $\nu : \setRp \rightarrow \setRp$ is regularly
  varying at $0$ if there is a real number $\beta \in \setR$ such that
  \begin{equation}
    \label{eq:rv_def}
    \forall y > 0, \quad \lim_{h \rightarrow 0^+}
    \nu(yh) / \nu(h) = y^{\beta}. 
  \end{equation}
  We denote by $\RV(\beta)$ the set of all such functions. A function
  in $\RV(0)$ is \emph{slowly varying}. 
\end{definition}

\begin{remark}
  Roughly speaking, a regularly varying function behaves as a power
  function times a slower term. Typical examples of such functions are
  $x^{\beta}$, $x^{\beta}(\log(1/x))^{\gamma}$ and more generally any
  power function times a $\log$ or compositions of $\log$ to some
  power. For other examples, see in the references. 
\end{remark}

\begin{definition}
  \label{def:rv_sigma_definition}
  If $\delta > 0$ and $\omega \in \RV(s)$ with $s > 0$ we define the
  class $\mc F_{\delta}(x_0, \omega)$ of all functions $f : \setR
  \raro \setR$ such that
  \begin{equation*}
    \forall h \leq \delta, \quad \inf_{P \in \mc P_{k}} \sup_{|x-x_0| \leq
      h} |f(x) - P(x-x_0)| \leq \omega(h),
  \end{equation*}
  where $k = \ppint{s}$ (the largest integer smaller than $s$) and
  $\mc P_k$ is the set of all the real polynomials with degree $k$. We
  define $\ellom(h) \eqdef \omega(h) h^{-s}$ the slow variation term
  of $\omega$. If $\alpha > 0$ we define
  \begin{equation*}
    \mc U(\alpha) \eqdef \bigl\{ f : [0,1] \raro \setR
    \text{ such that } \norminfty{f} \leq \alpha \bigr\}. 
  \end{equation*}
  Finally, we define
  \begin{equation*}
    \Sigma_{\delta, \alpha}(x_0, \omega) \eqdef \mc F_{\delta}(x_0,
    \omega) \cap \mc U(\alpha). 
  \end{equation*}
\end{definition}

\begin{remark}
  If $\omega(h) = r h^s$ for $r > 0$ we find back the classical
  H{\"o}lder regularity with radius $r$. In this sense, the class $\mc
  F_{\delta}(x_0, \omega)$ is a slight H{\"o}lder regularity
  generalisation.
\end{remark}

\subsection{Conditionally on the design}
\label{sec:conditional_on_the_design}

When nothing is known on the design density behaviour we can work
conditionally on the design. Let $\mf X_n$ be the sigma-algebra
generated by $X_1, \ldots, X_n$. We define
\begin{equation}
  \label{eq:H_n_def_adapt}
  H_{n,\omega} \eqdef \min \Bigl\{ h \in [0,1] \text{ such that }
  \omega(h) \geq \sigma\sqrt{N_{n,h}^{-1} \log n} \Bigr\},
\end{equation}
which is well defined for $n$ large enough (when $\omega(1) \geq
\sigma \sqrt{\log n / n}$). The quantity $H_{n,\omega}$ makes the
balance between the bias and the $\log$-penalised variance of $\wh
f_{h,\kpa}$ (see lemma \ref{lem:bias_variance_adapt}) and therefore
can be understood as the \emph{ideal adaptive bandwidth}, see Lepski
and Spokoiny (1997)\nocite{lepski_spok97} and Spokoiny
(1998)\nocite{spok98}.  The $\log$ term in \eqref{eq:H_n_def_adapt} is
the \emph{payment for adaptation}, see section
\ref{sec:payment_for_adaptation}. Let us define
\begin{equation*}
  H_{n, \omega}^{*} \eqdef \max \{ h \in \mc H | h \leq H_{n, \omega}
  \},
\end{equation*}
and
\begin{equation}
  \label{eq:R_n_omega_def}
  R_{n,\omega} \eqdef \sigma \sqrt{N_{n,H_{n,\omega}^*}^{-1} \log n}. 
\end{equation}
We define the diagonal matrix $\Lba_h \eqdef \text{diag}(
\normh{\phi_0}^{-1}, \ldots, \normh{\phi_{\kpa}}^{-1})$, the
symmetrical matrix $\mc G_h \eqdef \Lba_h \mbt X_h \Lba_h$ and
$\lba_{n, \omega} \eqdef \lba(\mc G_{H_{n,\omega}^*})$. We define the
event
\begin{equation}
  \label{eq:def_omega_X_invertible}
  \Omega_h \eqdef \{ X_1,\ldots,X_n \text{ are such that } \lba(\mb
  X_h) > N_{n,h}^{-1/2} \text{ and } N_{n,h} \geq 2 \}. 
\end{equation}
We note that $\Omega_h \in \mf X_n$ and $\mb X_h$ is invertible on
$\Omega_h$. The next result shows that, conditional on $\mf X_n$, $\wh
f_n(x_0) = \wh f_{\wh H_n, \kpa}(x_0)$ converges with the rate
$R_{n,\omega}$ simultaneously over any $\Sigma(x_0, \omega)$ when
$\omega \in \RV(s)$ with $0 < s \leq \kpa + 1$.

\begin{theorem}
  \label{thm:NA_adapt_upper_bound}
  If $\omega \in \RV(s)$, $0 < s \leq \kpa + 1$ and $\alpha > 0$ we
  have for any $n \geq \kpa +1$ on $\Omega_{H_{n,\omega}^*}$\tup:
  \begin{equation*}
    \sup_{\displaystyle f \in \Sigma_{H_{n, \omega}^*, \alpha}(x_0,
      \omega)} \Efm \bigl\{ R_{n, \omega}^{-p} | \wh f_n(x_0) - f(x_0)
    |^p | \mf X_n \bigr\} \leq c_1 \lba_{n, \omega}^{-p} + c_2
    (\alpha \vee 1)^p (\log n)^{-p / 2},
  \end{equation*}
  where $c_1 = c_1(p, \kpa, a)$ and $c_2 = c_2(p, \kpa, a, \sigma)$. 
\end{theorem}

We will see that the probability of the event $\Omega_{H_{n,
    \omega}^*}$ is large and that $\lba_{n, \omega}$ is positive with
a large probability when the design density is regularly varying (see
lemma \ref{lem:A_n_von_subset}). Note that the upper bound in theorem
\ref{thm:NA_adapt_upper_bound} is non asymptotic in the sense that it
holds for any $n \geq \kpa+1$. The random normalisation $R_{n,
  \omega}$ is similar to the one in Guerre (1999)\nocite{guerre99},
see section \ref{sec:comparison_with_existing_results} for more
details.

\subsection{Regularly varying design}
\label{sect:regularly_varying_design}

\begin{definition}
  \label{def:class_R_W_nu}
  For $\beta > -1$ and a neighbourhood $W$ of $x_0$ we define
  \begin{equation*}
    \mc R(x_0, \beta) \eqdef \bigl\{ \mu \text{ density such that }
    \exists \nu \in \RV(\beta) \, \forall x \in W, \quad \mu(x) =
    \nu(|x-x_0|) \bigr\}. 
  \end{equation*}
\end{definition}
We assume in all the following that $\mu \in \mc R(x_0, \beta)$ for
$\beta > -1$. Let $h_{n,\omega}$ be the smallest solution to
\begin{equation}
  \label{eq:bias_variance_det_adapt}
  \omega(h) = \sigma \sqrt{\frac{\log n}{2 n \int_0^h \nu(t) dt}},
\end{equation}
and 
\begin{equation}
  \label{eq:r_n_def}
  r_{n, \omega} \eqdef \omega(h_{n,\omega}).
\end{equation}
Equation \eqref{eq:bias_variance_det_adapt} can be viewed as the
deterministic counterpart to the equilibrium in
\eqref{eq:H_n_def_adapt}. We define $C_{\alpha, \beta} \eqdef (1 +
(-1)^{\alpha})\frac{\beta + 1}{\alpha + \beta+1}$ and the matrix $\mc
G$ with entries $(\mc G)_{j, l} \eqdef \frac{C_{j + l,
    \beta}}{\sqrt{C_{2j, \beta} C_{2l, \beta}}}$ for $0 \leq j,l \leq
\kpa$ and $\lba_{\kpa,\beta} \eqdef \lba(\mc G)$. It is easy to see
that $\lba_{\kpa, \beta} > 0$.

\begin{theorem}
  \label{thm:asympt_adaptive_upper_bound}
  If
  \begin{itemize}
  \item $\kpa \in \setN$, $\beta > -1$, $\alpha > 0$ and $\varrho >
    1$,
  \item $\omega \in \RV(s)$ for $0 < s \leq \kpa + 1$,
  \end{itemize}
  then the estimator $\wh f_n(x_0) = \wh f_{\kpa, \wh H_n}(x_0)$ with
  the grid $\mc H = \mc H_1^{\tup{arith}}$ satisfies
  \begin{equation}
    \label{eq:adapt_upper_bound_asympt}
    \forall \mu \in \mc R(x_0, \beta) \quad \limsup_n
    \sup_{\displaystyle f \in \Sigma_{\varrho h_{n, \omega},
        \alpha}(x_0, \omega)} \Efm \bigl\{ r_{n, \omega}^{-p}
    |\wh f_n(x_0) - f(x_0)|^p \bigr\} \leq C \lba_{\kpa, \beta}^{-p},
  \end{equation}
  where $C = C(p, \kpa)$. Moreover, we have
  \begin{equation}
    \label{eq:upper_bound_r_n_equiv_adapt}
    r_{n, \omega} \sim \sigma^{2s/(1+2s+\beta)} (\log n/n)^{s / (1 +
      2s + \beta)} \ell_{\omega,\nu}(\log n / n) \text{ as } n \raro
    +\infty,
  \end{equation}
  where $\ell_{\omega,\nu}$ is slowly varying. 
\end{theorem}

\begin{remark}
  When $\omega(h) = r h^s$ (H{\"o}lder regularity) we have more
  precisely
  \begin{equation*}
    r_{n, \omega} \sim \sigma^{2s / (1 + 2s + \beta)} r^{(1 +
      \beta) / (1 + 2s + \beta)} (\log n/n)^{s/(1 + 2s + \beta)}
    \ell_{s,\nu}( \log n / n) \text{ as } n \raro +\infty. 
  \end{equation*}
  Note that $\ell_1(h) = \ell_{\omega, \nu}(h \log(1 / h))$ is also
  slowly varying, thus $\ell_{1}(1 / n) = \ell_{\omega, \nu} (\log n /
  n)$ is a slow term. 
\end{remark}

\subsection{Convergence rates examples}
\label{sec:example_of_rate_of_convergence}

Let $\beta > -1$, $r, s$ be positive and $\alpha, \gamma$ be any real
numbers. If we take $\nu$ such that $\int_0^h \nu(t) dt = h^{\beta +
  1}(\log(1 / h))^{\alpha}$ and $\omega(h) = r h^s
(\log(1/h))^{\gamma}$ then we find that (see section
\ref{sec:computation_exemple} for the computation details)
\begin{equation}
  \label{eq:convergence_rate_exemple}
  r_{n, \omega} \sim \sigma^{2s/(1+2s+\beta)} r^{(\beta + 1)/(1 + 2s +
    \beta)} \bigl(n (\log n)^{ \alpha - 1 - \gamma (1 + \beta) / s}
  \bigr)^{-s/(1 + 2s + \beta)},
\end{equation}
where $a_n \sim b_n$ mean $\lim_{n \raro +\infty} a_n / b_n = 1$. This
rate has to be compared with the minimax rate from Ga\"iffas
(2004)\nocite{gaiffas04a}:
\begin{equation*}
  \sigma^{2s/(1+2s+\beta)} r^{(\beta + 1)/(1 + 2s + \beta)} \bigl(n
  (\log n)^{ \alpha - \gamma (1 + \beta) / s} \bigr)^{-s/(1 + 2s
    + \beta)},
\end{equation*}
where the only difference is the $\alpha$ instead of $\alpha -1$ in
the $\log$ exponent. This loss is the \emph{payment for adaptation}
and is unavoidable in view of theorem \ref{thm:lower_bound} and its
corollary. See section \ref{sec:optimality} for more details. 

In the classical case, namely when the design is non-degenerate and
$f$ is H{\"o}lder ($\omega(h) = r h^s$ and $\alpha = \beta = \gamma=
0$) we find the usual pointwise minimax adaptive rate (see Lepski
(1990)\nocite{lepski90}, Brown and Low
(1996)\nocite{brown_low_constrained96}):
\begin{equation*}
  \sigma^{2s / (1 + 2s)} r^{1 / (1 + 2s)} (\log n / n)^{s/(1+2s)}. 
\end{equation*}
When the design is again non-degenerate and the continuity modulus is
equal to $\omega(h) = r h^s (\log(1/h))^{-s}$ we find a convergence
rate equal to
\begin{equation*}
  \sigma^{2s / (1 + 2s)} r^{1 / (1 + 2s)} n^{-s/(1+2s)},
\end{equation*}
which is the usual minimax rate, without the $\log$ term for payment
for adaptation. Actually, this is a "toy" example since we have asked
for more regularity than in the H{\"o}lder regularity. Note that in
the degenerate design case, when $\alpha$ and $\gamma$ are such that
$\alpha = 1 + \gamma (1+\beta)/ s$ there is again no extra $\log$
factor.

\section{Optimality}
\label{sec:optimality}

\subsection{Payment for adaptation}
\label{sec:payment_for_adaptation}

The convergence rate of a linear estimator with an adaptive bandwidth
choice can be well explained with a balance equation between its bias
and variance terms. In our context this equation is
\begin{equation*}
  \omega(h) = \frac{\sigma}{\sqrt{N_{n,h}}},
\end{equation*}
(see lemma \ref{lem:bias_variance_adapt}) and a deterministic
counterpart of this equilibrium is
\begin{equation}
  \label{eq:bias_variance_det_no_log}
  \omega(h) = \frac{\sigma}{\sqrt{2 n \int_0^h \nu(t) dt}},
\end{equation}
see lemma \ref{lem:Nh_prob_equiv}. We proved in Ga\"iffas
(2004)\nocite{gaiffas04a} that the minimax rate $\psi_{n, \omega}$
over $\Sigma_{\delta, \alpha}(x_0, \omega)$ is given by
\begin{equation}
  \label{eq:psi_n_def}
  \psi_{n, \omega} =  \omega( \gamma_{n, \omega}),
\end{equation}
where $\gamma_{n, \omega}$ is the smallest solution to
(\ref{eq:bias_variance_det_no_log}). In a model with
\emph{homogeneous} information (the white noise or the regression
model with an equidistant design) we know that such a balance equation
cannot be realized: an adaptive estimator to the unknown smoothness
without loss of efficiency is not possible for pointwise estimation,
even if we know that the function belongs to one of two H{\"o}lder
classes, see Lepski (1990)\nocite{lepski90}, Brown and Low
(1996)\nocite{brown_low_constrained96} and Lepski and Spokoiny (1997)
\nocite{lepski_spok97}.  This means that local adaptation cannot be
achieved for \emph{free}: we have to pay an extra $\log$ factor in the
convergence rate, at least of order $(\log n)^{2s/(1+2s)}$ when
estimating a H{\"o}lder function with smoothness $s$. The authors call
this phenomenon \emph{payment for adaptation}.  We intend here to
generalise this result to the regression with a degenerate random
design.

\subsection{Superefficiency}
\label{sec:a_lower_bound}

Let $s$, $r' < r$, be positive and $\delta \leq 1$, $p > 1$. We take
$\omega(h) = r h^s$, $\omega'(h) = r' h^s$ and the minimax rate
$\psi_{n, \omega}$ defined by (\ref{eq:psi_n_def}). In view of lemma
\ref{lem:rv_r_n_equiv_adapt} we have
\begin{equation}
  \label{eq:psi_n_equiv_lower_bound_part}
  \psi_{n, \omega} \sim c_{s, \beta, \sigma, r} 
  n^{ -s / (1 + 2 s + \beta)} \ell_{s, \nu}(1 / n)
  \text{ as } n \raro +\infty.
\end{equation}
We recall that in view of theorem
\ref{thm:asympt_adaptive_upper_bound} the "adaptive" rate $r_{n,
  \omega}$ defined by (\ref{eq:r_n_def}) is attained by the adaptive
procedure $\wh f_n(x_0)$ simultaneously over several classes
$\Sigma_{\delta, \alpha}(x_0, \omega)$ with $\omega \in \RV(s)$ for
any regularity $s \in (0, \kpa+1]$ and that
\begin{equation}
  \label{eq:r_n_equiv_lower_bound_part}
  r_{n, \omega} \sim c_{s, \beta, \sigma, r} 
  (\log n/n)^{s / (1 + 2s + \beta)} \ell_{s, \nu}(\log
  n / n) \text{ as } n \raro +\infty. 
\end{equation}

\begin{theorem}
  \label{thm:lower_bound}
  If an estimator $\wh f_n$ based on \eqref{eq:regression_model} is
  asymptotically minimax over $\mc F_{\delta}(x_0, \omega)$, that is
  \begin{equation*}
    \limsup_n \sup_{\displaystyle f \in \mc F_{\delta}(x_0, \omega)}
    \psi_{n, \omega}^{-p} \,\, \Efm \{ | \wh f_n(x_0) - f(x_0) |^p \} <
    +\infty,
  \end{equation*}
  and if this estimator is superefficient at a function $f_0 \in \mc
  F_{\delta}(x_0, \omega')$ in the sense that there is $\gamma > 0$
  such that
  \begin{equation}
    \label{eq:the_estimator_is_superefficient}
    \limsup_n \psi_{n, \omega}^{-p} \, n^{\gamma p} \,\, \mbb E_{f_0, \mu}^n
    \{ | \wh f_n(x_0) - f_0(x_0) |^p \} < +\infty,
  \end{equation}
  then we can find a function $f_1 \in \mc F_{\delta}(x_0, \omega)$
  such that
  \begin{equation*}
    \liminf_n r_{n, \omega}^{-p} \,\, \mbb E_{f_1, \mu}^n \{ | \wh f_n(x_0)
    - f_1(x_0) |^p \} > 0. 
  \end{equation*}
\end{theorem}

This theorem is a generalisation of a result by Brown and Low
(1996)\nocite{brown_low_constrained96} for the degenerate random
design case. Of course, when the design is non-degenerate ($0 <
\mu(x_0) < +\infty$) the theorem remains valid and the result is
barely the same as in Brown and Low
(1996)\nocite{brown_low_constrained96} with the same rates.

The theorem \ref{thm:lower_bound} is a lower bound for a
superefficient estimator. Actually, the most interesting result for
our problem is the next corollary.

\subsection{An adaptive lower bound}
\label{sec:adaptive_lower_bound}

Let $0 < r_2 < r_1 < +\infty$ and $0 < s_1 < s_2 < +\infty$ be such
that $\ppint{s_1} = \ppint{s_2} = k$. If $\omega_i(h) = r_i h^{s_i}$
we denote $\mc F_i \eqdef \mc F_{\delta}(x_0, \omega_i)$. Let
$\psi_{n, i}$ be the minimax rate defined by (\ref{eq:psi_n_def}) over
$\mc F_i$ for $i=1, 2$ and $r_{n,1}$ be defined by (\ref{eq:r_n_def})
with $\omega = \omega_1$ (the "adaptive" rate when the class is $\mc
F_1$). Note that $\psi_{n, i}$ satisfies
\eqref{eq:psi_n_equiv_lower_bound_part} with $s = s_i$ and $r_{n, 1}$
satisfies \eqref{eq:r_n_equiv_lower_bound_part} with $s = s_1$.

\begin{corollary}
  \label{corollary_lower_bound}
  If an estimator $\wh f_n$ is asymptotically minimax over $\mc F_1$
  and $\mc F_2$, that is for $i= 1, 2$\tup:
  \begin{equation}
    \label{eq:minimax__estimator_over_F_1_and_F_2}
    \limsup_n \sup_{\displaystyle f \in \mc F_i} \psi_{n,i}^{-p} \,\,
    \Efm\{ |\wh f_n(x_0) - f(x_0) |^p \} < +\infty,
  \end{equation}
  then this estimator also satisfies
  \begin{equation}
    \label{eq:corollary_lower_bound}
    \liminf_n \sup_{\displaystyle f \in \mc F_1} r_{n,1}^{-p} \,\, \Efm \big\{ |\wh
    f_n(x_0) - f(x_0) |^p \big\} > 0.
  \end{equation}
\end{corollary}
Note that (\ref{eq:corollary_lower_bound}) contradicts
(\ref{eq:minimax__estimator_over_F_1_and_F_2}) for $i = 1$ since
$\lim_n \psi_{n, 1} / r_{n,1} = 0$, thus \emph{\textbf{there is no
    pointwise minimax adaptive estimator over two such classes $\mc
    F_1$ and $\mc F_2$ and the best achievable rate is $r_{n, i}$.}}
The corollary \ref{corollary_lower_bound} is an immediate consequence
of theorem \ref{thm:lower_bound}. Clearly, $\mc F_2 \subset \mc F_1$
thus equation \eqref{eq:minimax__estimator_over_F_1_and_F_2} entails
that $\wh f_n$ is superefficient at any function $f_0 \in \mc F_2$.
More precisely, $\wh f_n$ satisfies
\eqref{eq:the_estimator_is_superefficient} with $\gamma = \frac{(s_2 -
  s_1)(\beta + 1)}{2 (1 + 2 s_1 + \beta)(1 + 2 s_2 + \beta)} > 0$
since $n^{- \gamma} \ell(1 / n) \raro 0$ where $\ell \eqdef \ell_{s_1,
  \nu} / \ell_{s_2, \nu}$ and $\ell \in \RV(0)$.


\section{Simulations}
\label{sec:simulations}

\subsection{Implementation of the procedure}
\label{sec:implementatio_of_the_procedure}

For the estimation at a point $x$, the procedure
\eqref{eq:H_n_hat_def} selects the best symmetrical interval $I = [x-
h, x+h]$ among several $h$ in the grid $\mc H$. We have implemented
this procedure with non symmetrical intervals, which is a procedure
similar to the one in Spokoiny (1998)\nocite{spok98}. First, we define
similarly to section \ref{sec:local_polynomial_estimation} for any $I
\subset [0, 1]$ the scalar product
\begin{equation*}
  \prodsca{f}{g}_I \eqdef \sum_{X_i \in I} f(X_i) g(X_i),
\end{equation*}
(it is convenient in this part to remove the normalisation term $N_{n,
  h}$ from the definition of the scalar product) and similarly to
\eqref{eq:matrix_X_vector_Y_def} we define the matrix $\mb X_I$ with
entries $(\mb X_I)_{j, l} = {\prodsca{\phi_j}{\phi_l}}_I$ for $0 \leq
j, l \leq \kpa$. We define in the same way $\mb Y_I$, and $\wh \tta_I$
is defined as the solution to
\begin{equation*}
  \mb X_I \tta = \mb Y_I.
\end{equation*}
Note that if $J \subset [0, 1]$, the vector $\mb F_{I, J}$ with
coordinates
\begin{equation*}
  (\mb F_{I, J})_j = {\prodsca{\wh f_{I, \kpa} - \wh f_{J,
        \kpa}}{\phi_j}}_{J} / \norm{\phi_j}_J
\end{equation*}
(for $0 \leq j \leq \kpa$) satisfies
\begin{equation*}
  \mb F_{I, J} = \mb H_J (\wh \tta_I - \wh \tta_J),
\end{equation*}
where $\mb H_J$ is defined as the matrix with entries for $0 \leq j, l
\leq \kpa$
\begin{equation*}
  (\mb H_J)_{j,l} \eqdef \frac{\sum_{ X_i \in J} (X_i - x)^{j
      + l}} {\sqrt{\sum_{X_i \in J} (X_i - x)^{2j}}}. 
\end{equation*}
The main steps of the procedure for the estimation at a point $x$ are
then:
\begin{enumerate}
\item choose parameters $a > 1$, $\kpa \in \setN$ and $m \geq \kpa +
  1$,
\item sort the $(X_i, Y_i)$ in $(X_{(i)}, Y_{(i)})$ such that $X_{(i)}
  \leq X_{(i+1)}$,
\item find $j$ such that $x \in [ X_{(j)}, X_{(j+1)}]$ and $\#\{ X_i |
  X_i \in [ X_{(j)}, X_{(j+1)}] \} = m$,
\item build
  \begin{equation*}
    \mc G^- = \bigcup_{i=0}^{ [\log_a(j + 1)] } \{ X_{( j + 1 - [a^i]
      )} \}, \hspace{1cm} \mc G^+ = \bigcup_{i=0}^{ [\log_a(n-j)] }
    \{ X_{( j + [a^i] )} \},
  \end{equation*}
\item compute $\wh \tta_I$ and $\mb H_I$ for all $I \in \mc G \eqdef
  \mc G^- \times \mc G^+$,
\item if $N_{n, I} \eqdef \#\{ X_i | X_i \in I \}$, find
  \begin{equation*}
    \wh I \eqdef \argmax_{I \in \mc G} \big\{ N_{n, I} \text{ such
      that } \forall J \subset I, J \in \mc G, \;
    \norminfty{\mb H_J (\wh \tta_I - \wh \tta_J) } \leq T_{I, J}
    \big\},
  \end{equation*}
  where $\norminfty{\cdot}$ stands for the sup norm in
  $\setR^{\kpa+1}$ and
  \begin{equation*}
    T_{I, J} = \wh \sigma (1 + \sqrt{\kpa+1}) \sqrt{\log N_{n, I}} +
    \sqrt{(1 + a)} \sqrt{ (N_{n, J} / N_{n, I}) \log n },
  \end{equation*}
  with $\wh \sigma$ for instance given by \eqref{eq:sigma_estim}.
\item return the first coordinate of $\wh \tta_{\wh I}$.
\end{enumerate}
This procedure uses a geometrical grid, thus it is computationally
feasible for reasonable choices of $a$ ($a=1.05$ is used in the next
section). The main steps of the procedure with an arithmetical grid
are the same with a modification of the threshold, see
\eqref{eq:T_n_def}. The procedure is implemented in \texttt{C++} and
is quite fast: it takes few seconds to recover the whole function at
$300$ points on a modern computer.


\subsection{Numerical illustrations}
\label{sec:subsect_simulations}

\begin{figure}[htbp]
  \begin{center}
    \includegraphics[width =
    \figurelength]{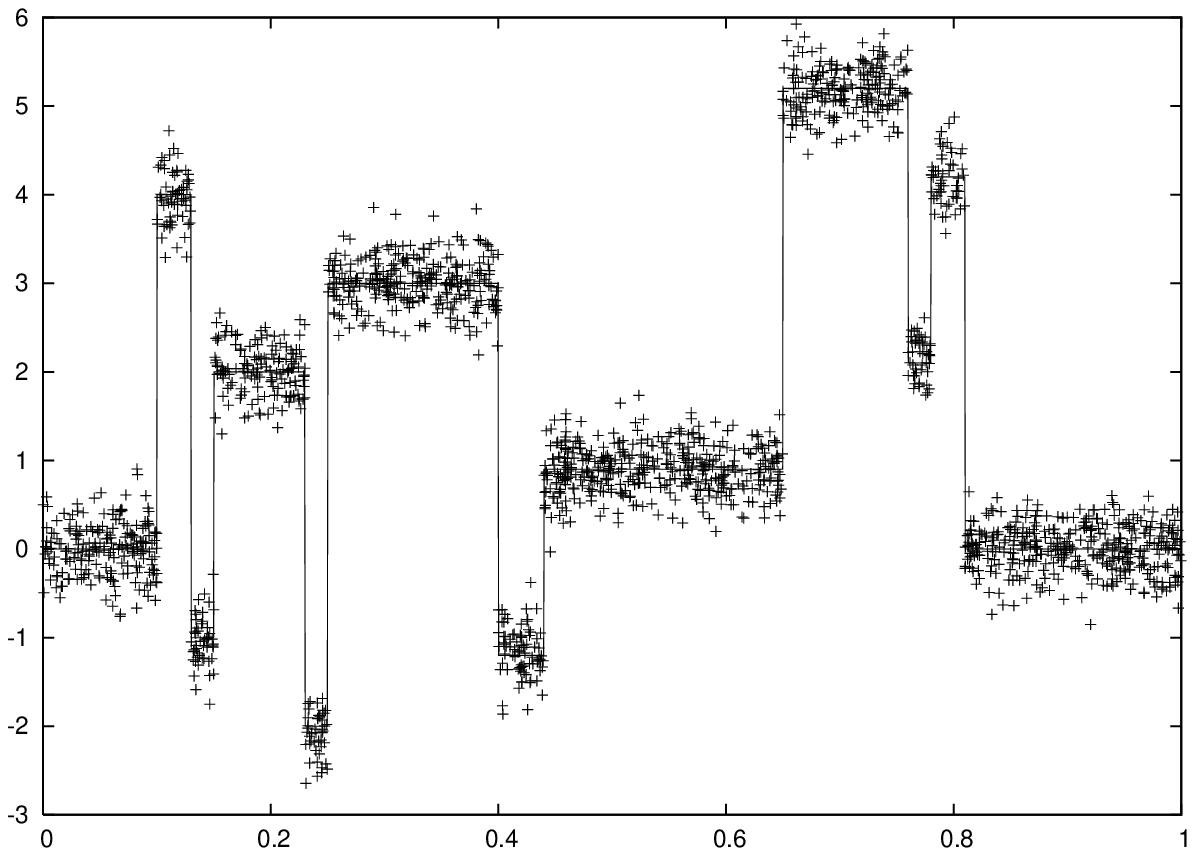}
    \includegraphics[width =
    \figurelength]{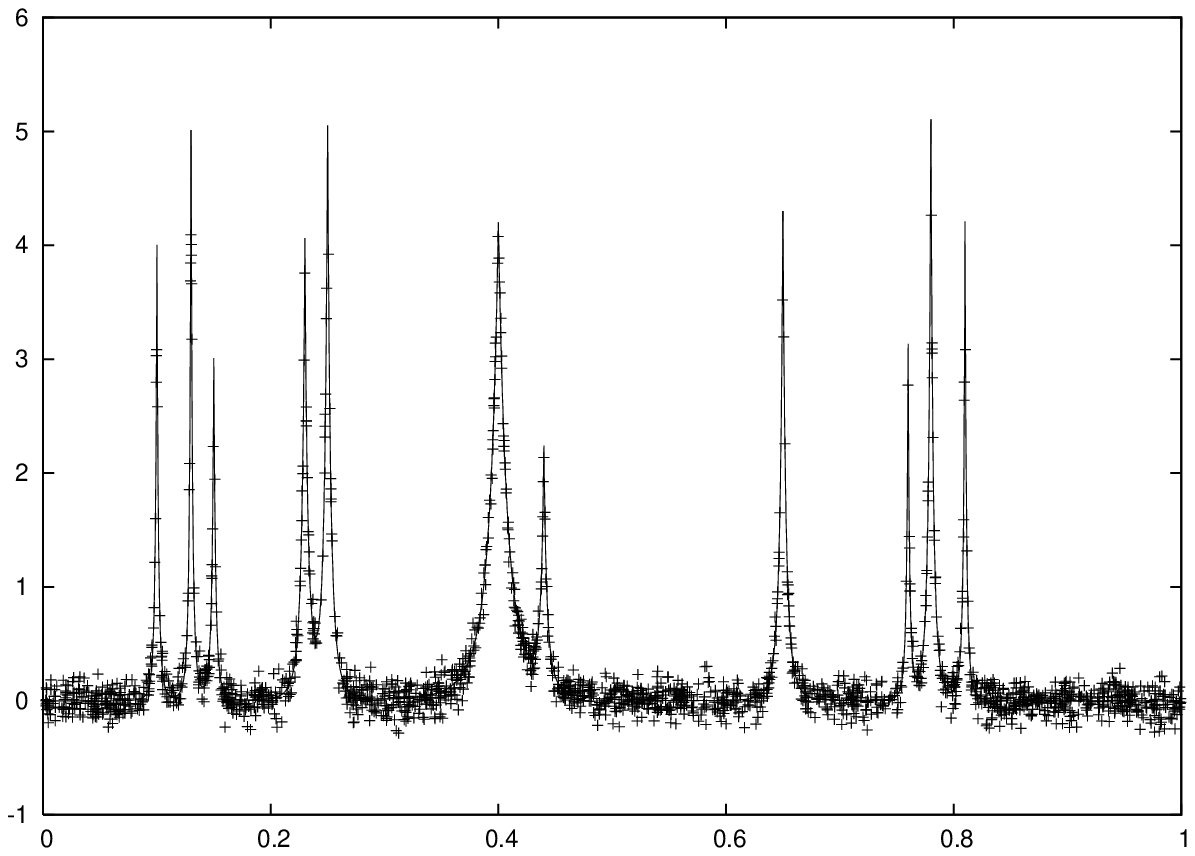}
    \includegraphics[width =
    \figurelength]{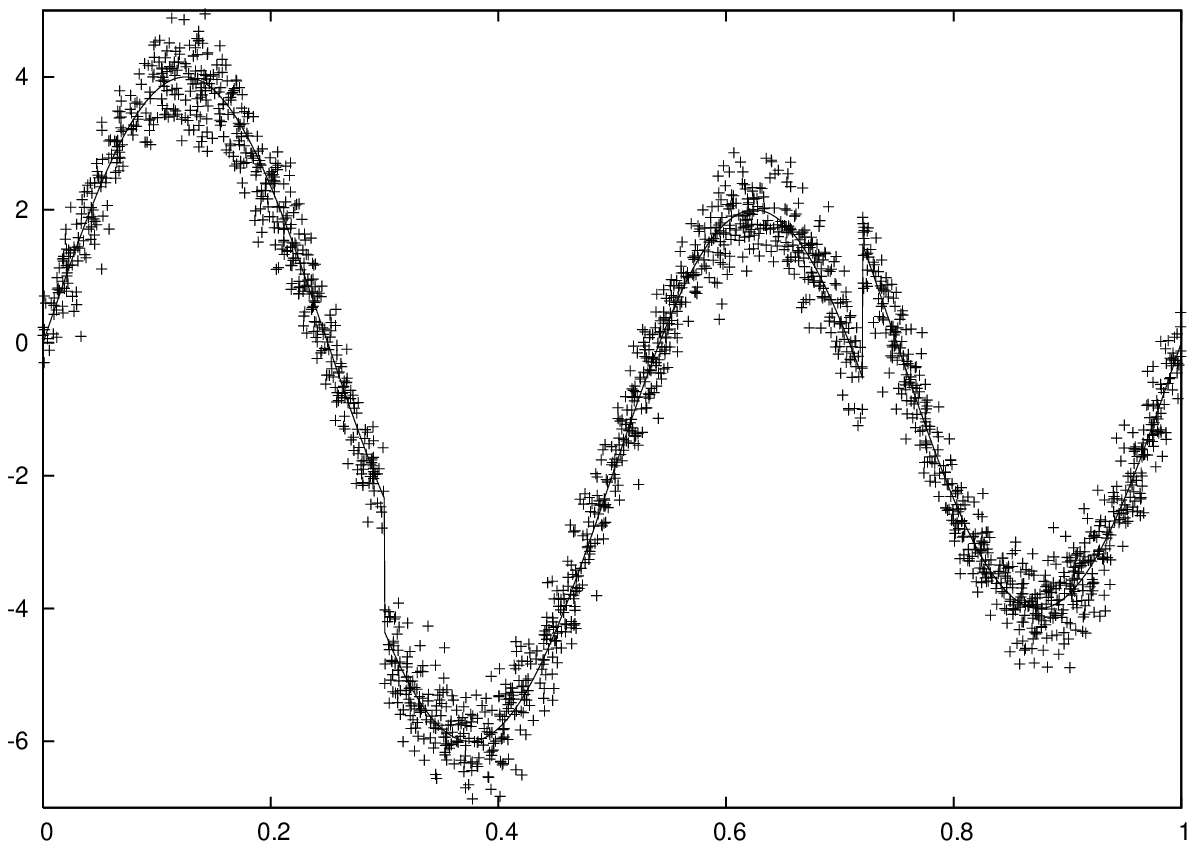}
    \includegraphics[width =
    \figurelength]{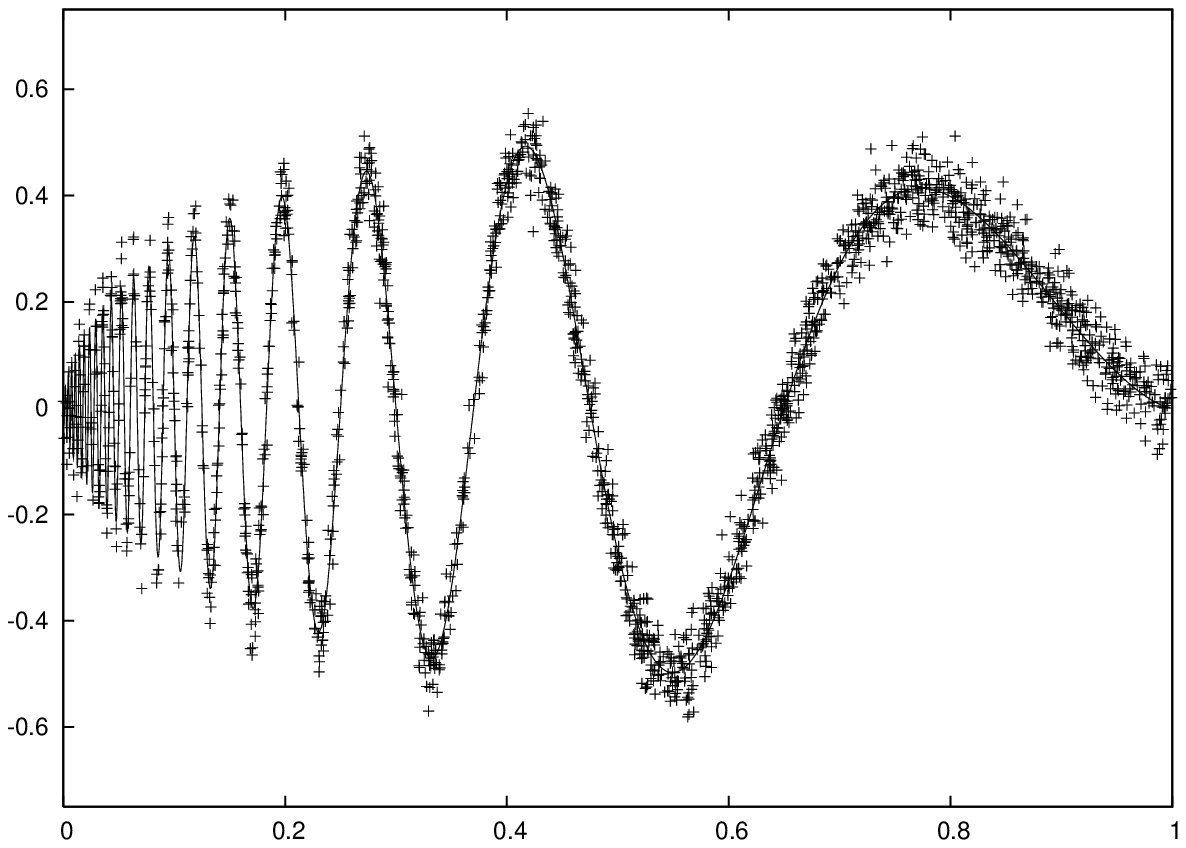}
  \end{center}
  \caption{Blocks, bumps, heavysine and doppler with Gaussian noise
    and uniform design.}
  \label{fig_datasets_unif_design}
\end{figure}

\begin{figure}[htbp]
  \begin{center}
    \includegraphics[width =
    \figurelength]{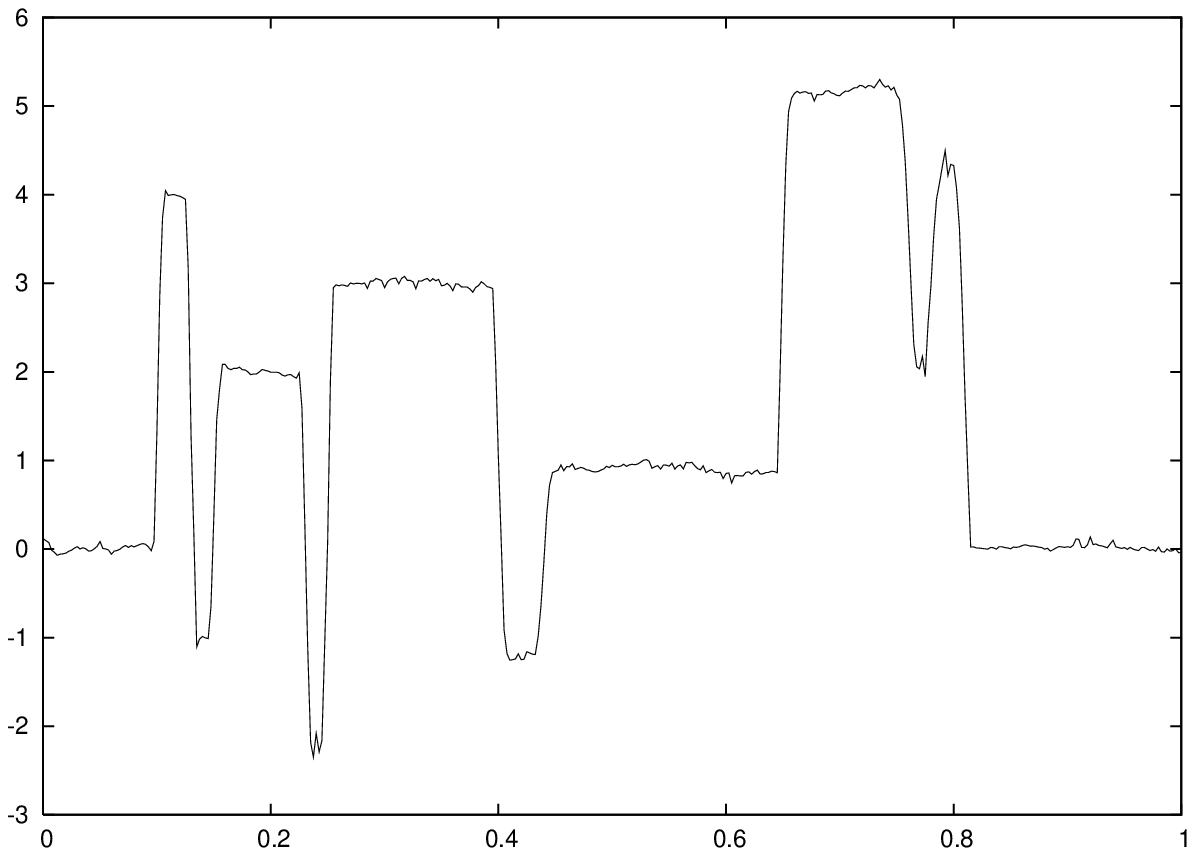}
    \includegraphics[width =
    \figurelength]{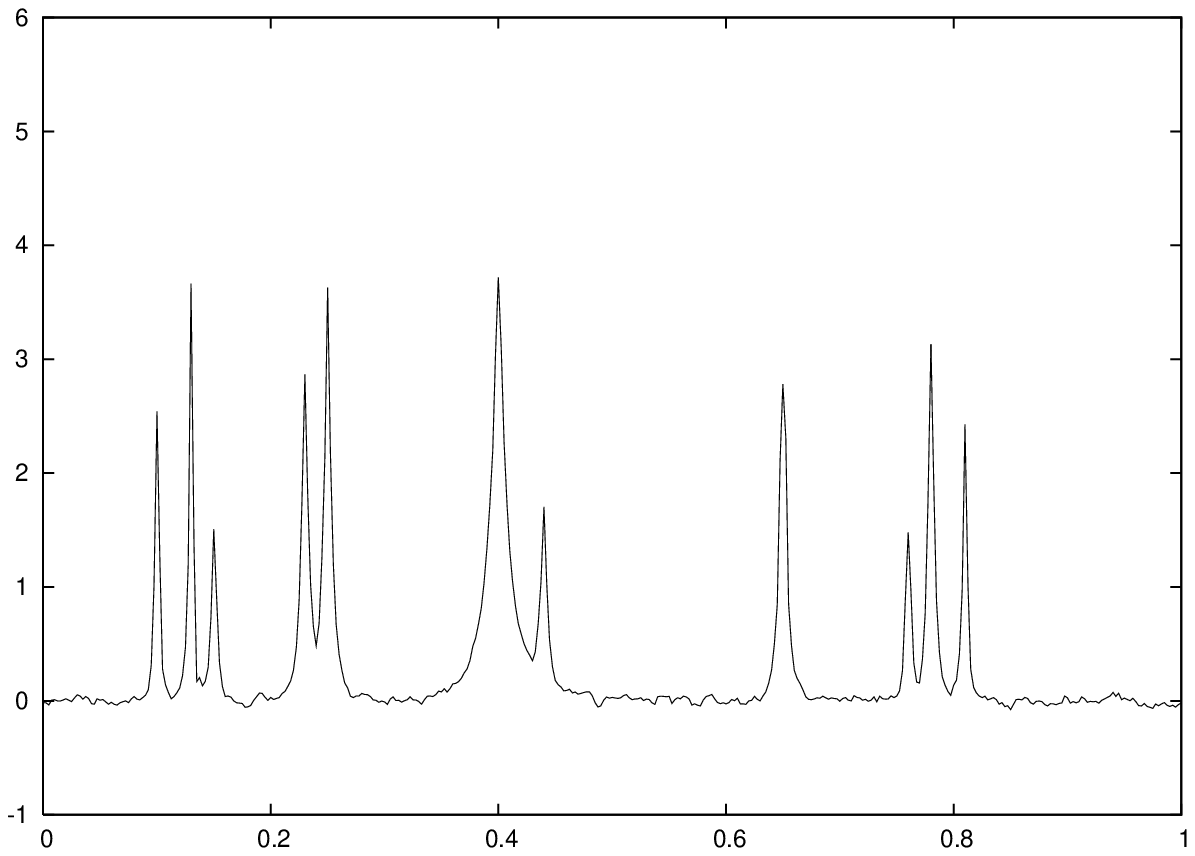}
    \includegraphics[width =
    \figurelength]{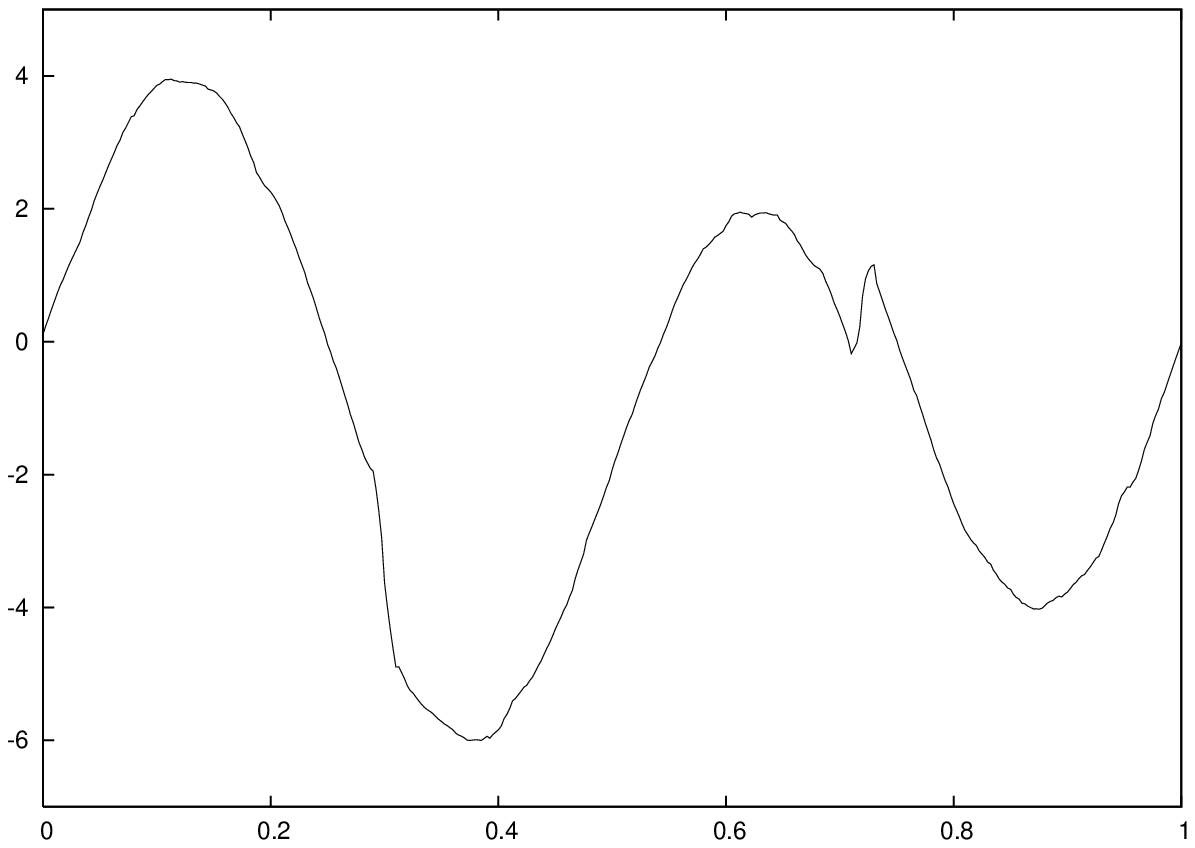}
    \includegraphics[width =
    \figurelength]{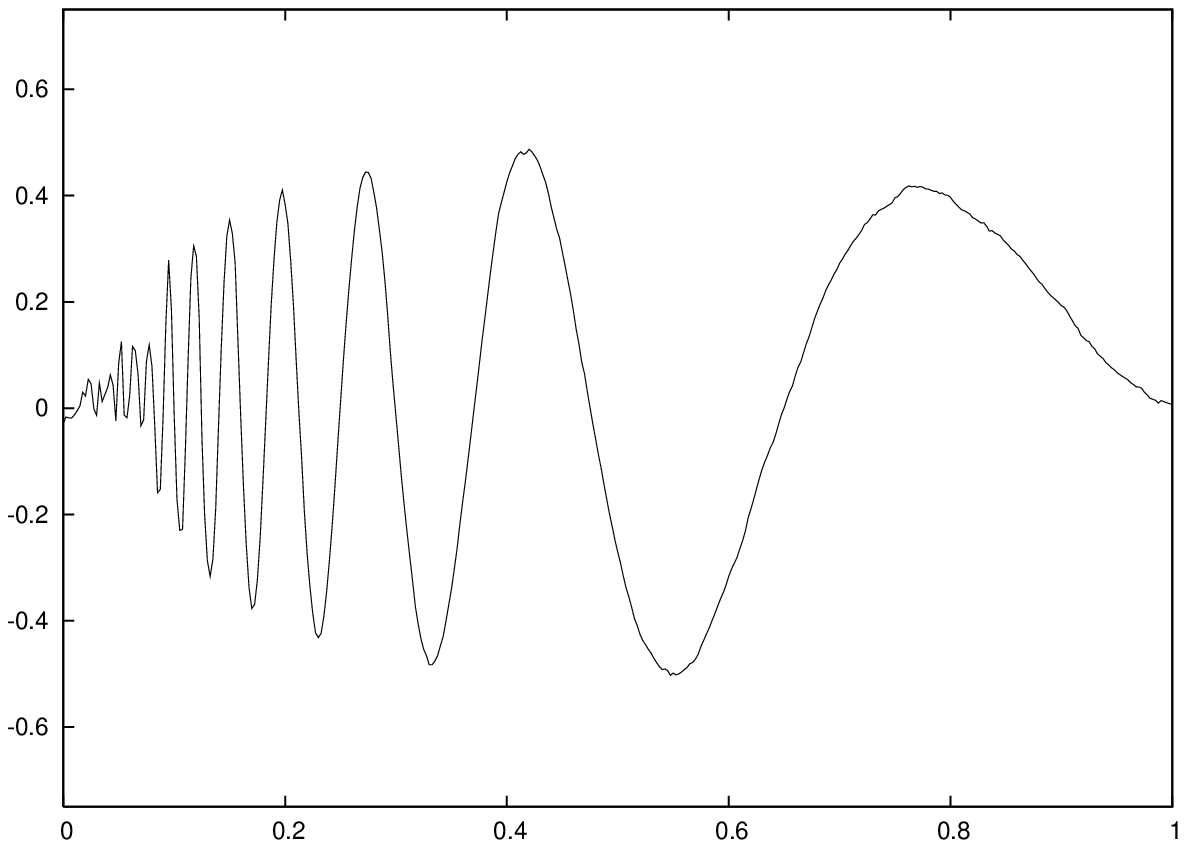}
  \end{center}
  \caption{Estimates based on the datasets in figure
    \ref{fig_datasets_unif_design}.}
  \label{fig_estimates_unif_design}
\end{figure}

We use for our simulations the target functions from Donoho and
Johnstone (1994)\nocite{donoho_johnstone94}. These functions are
commonly used as benchmarks for adaptive estimators. We show in figure
\ref{fig_datasets_unif_design} the target functions and datasets with
a uniform random design. The noise is Gaussian with $\sigma$ chosen to
have (root) signal-to-noise ratio $7$. The sample size is $n = 2000$.
We show the estimates in figure \ref{fig_estimates_unif_design}. For
all estimates we take $\kpa = 2$, $a = 1.05$ and $m = 25$. We estimate
at each point $x = j/300$ with $j = 0,\ldots,300$.


\newpage

\begin{figure}[htbp]
  \begin{center}
    \includegraphics[width =
    \figurelength]{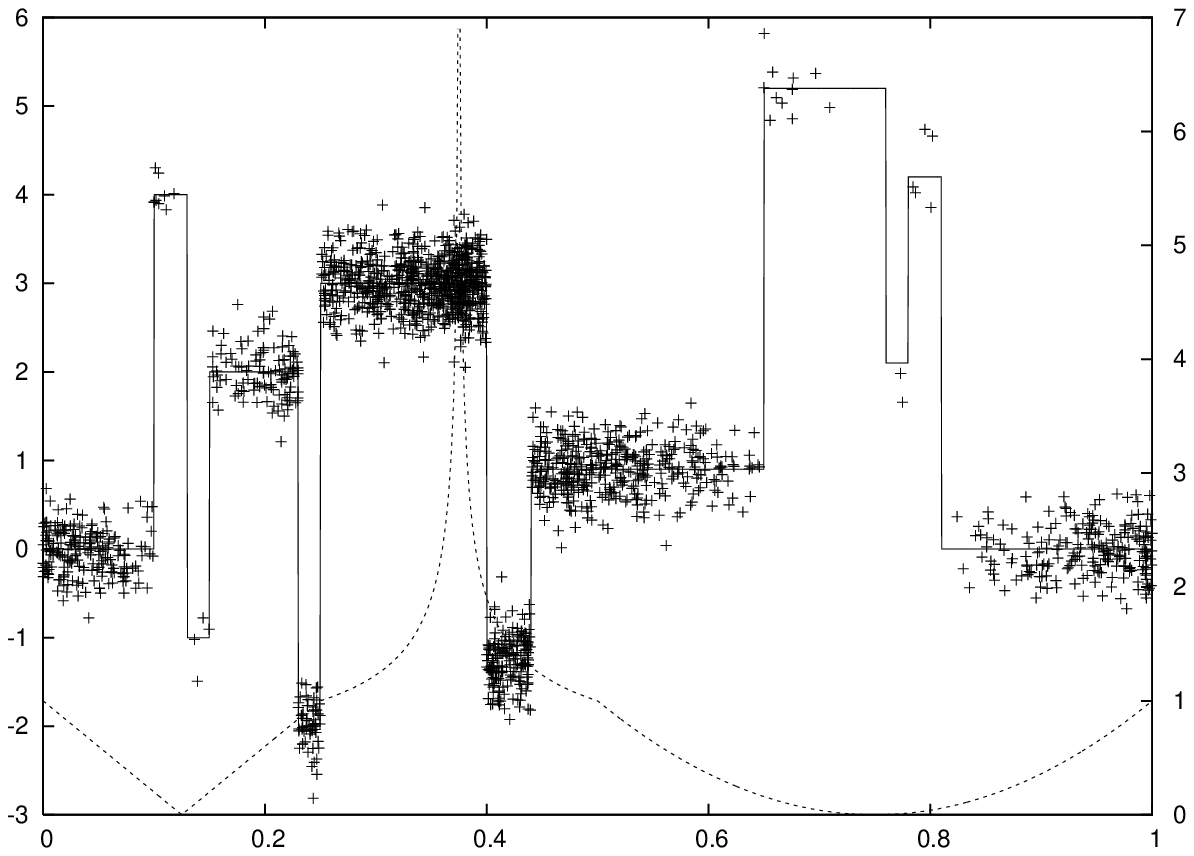}
    \includegraphics[width =
    \figurelength]{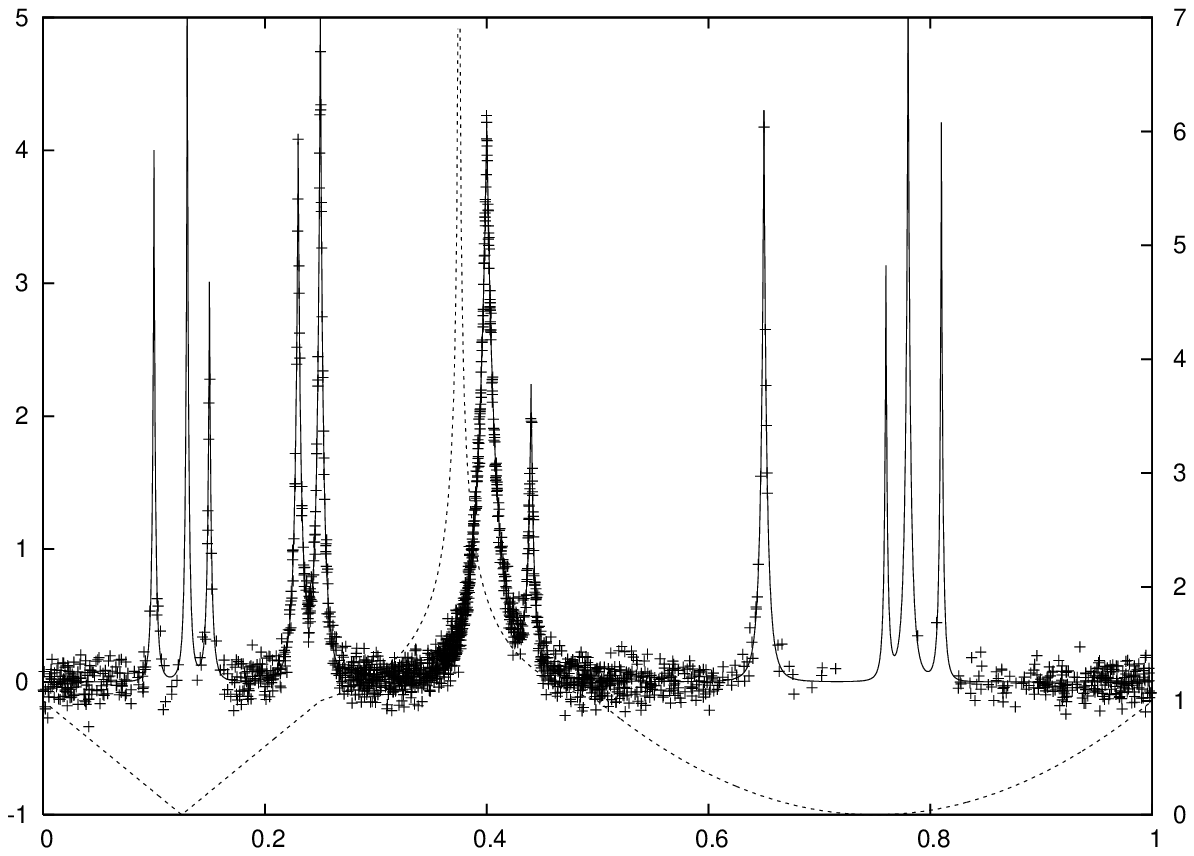}
    \includegraphics[width =
    \figurelength]{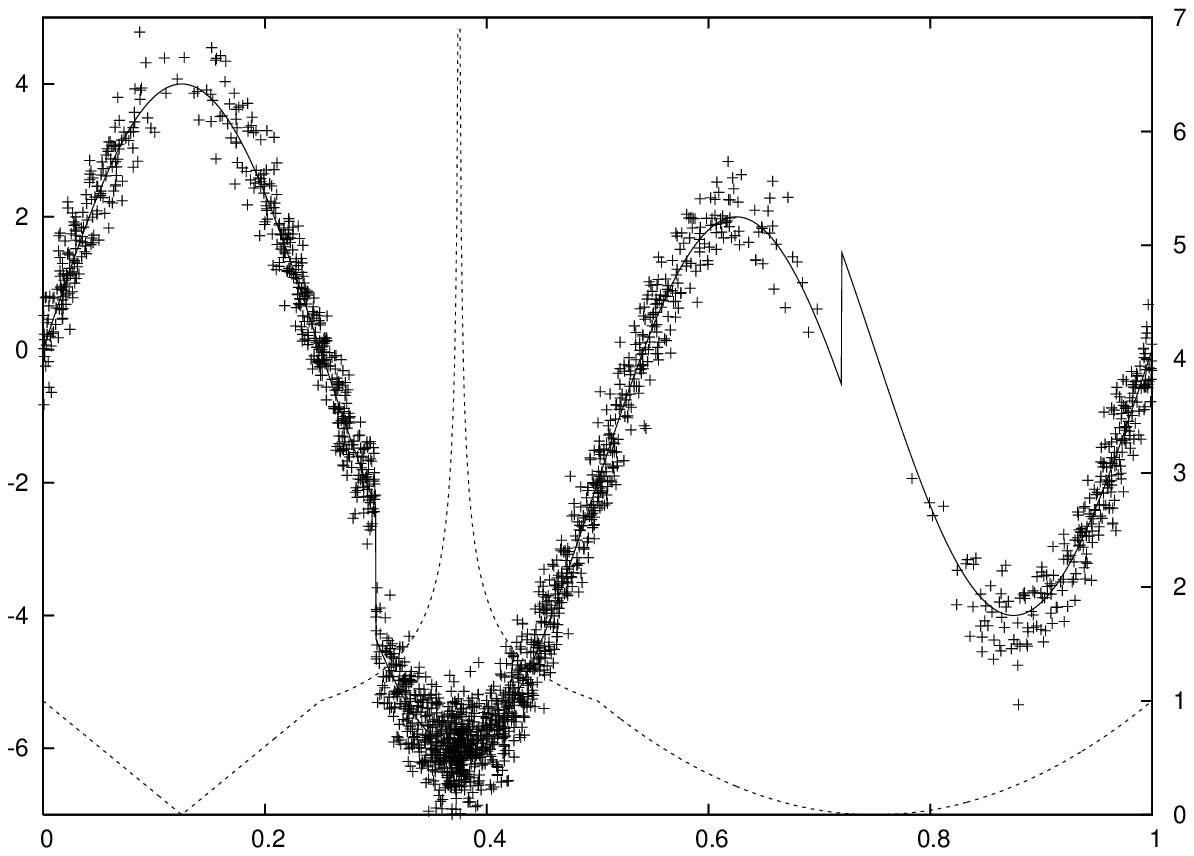}
    \includegraphics[width =
    \figurelength]{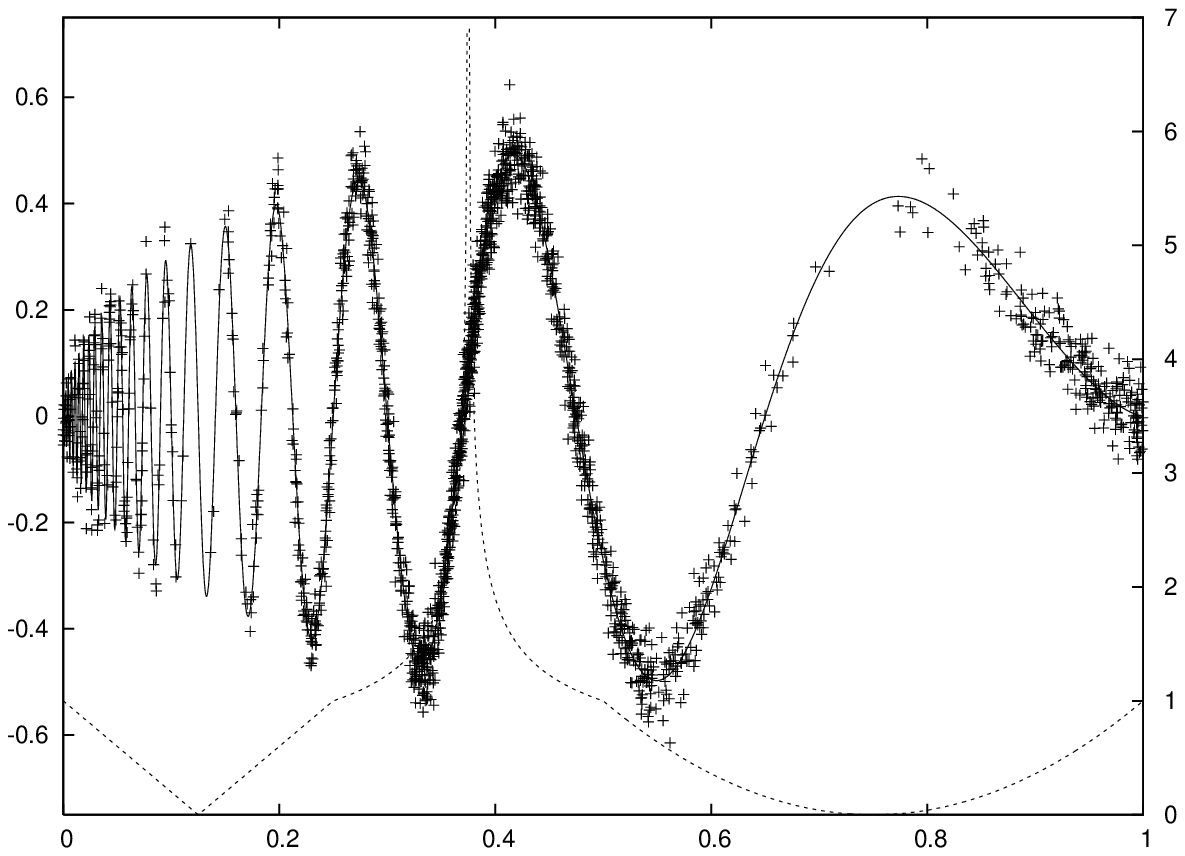}
  \end{center}
  \caption{Blocks, bumps, heavysine and doppler with Gaussian noise
    and non-uniform design.}
  \label{fig_datasets_example1_design}
\end{figure}

\begin{figure}[htbp]
  \begin{center}
    \includegraphics[width =
    \figurelength]{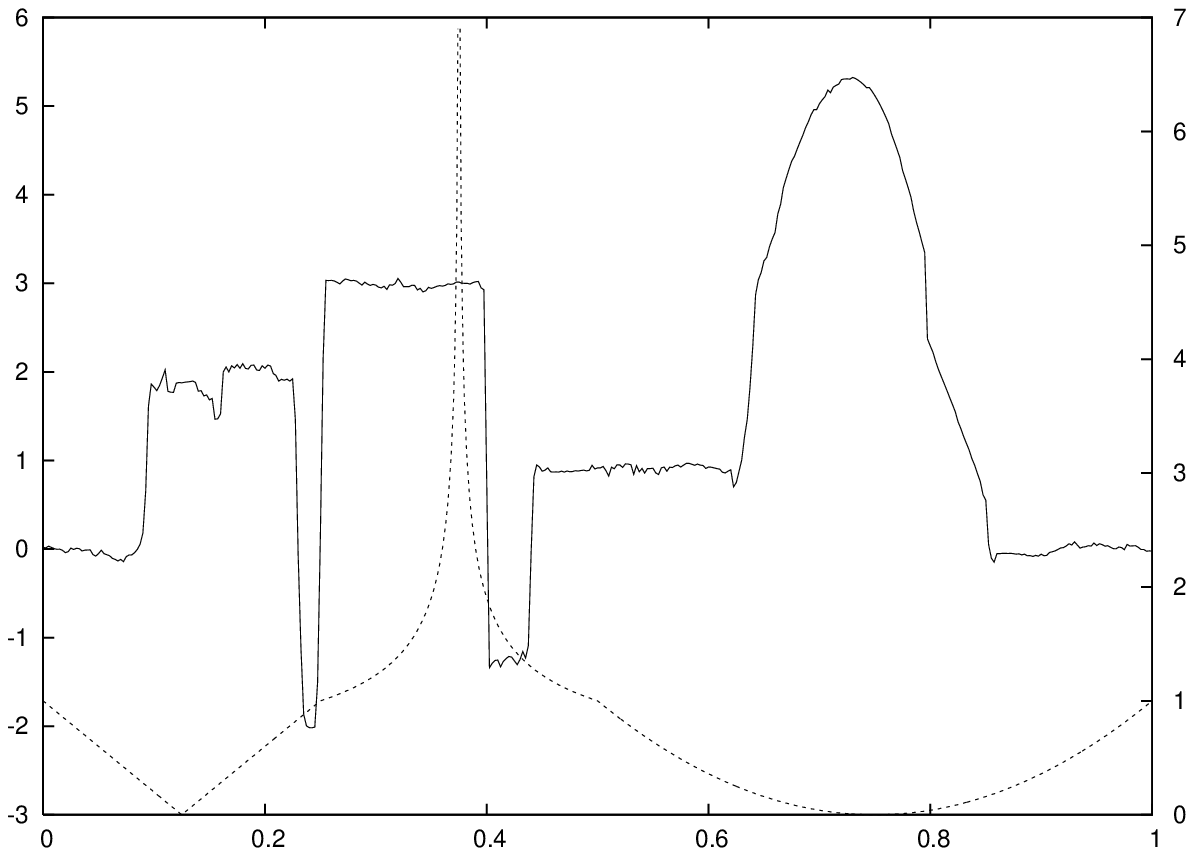}
    \includegraphics[width =
    \figurelength]{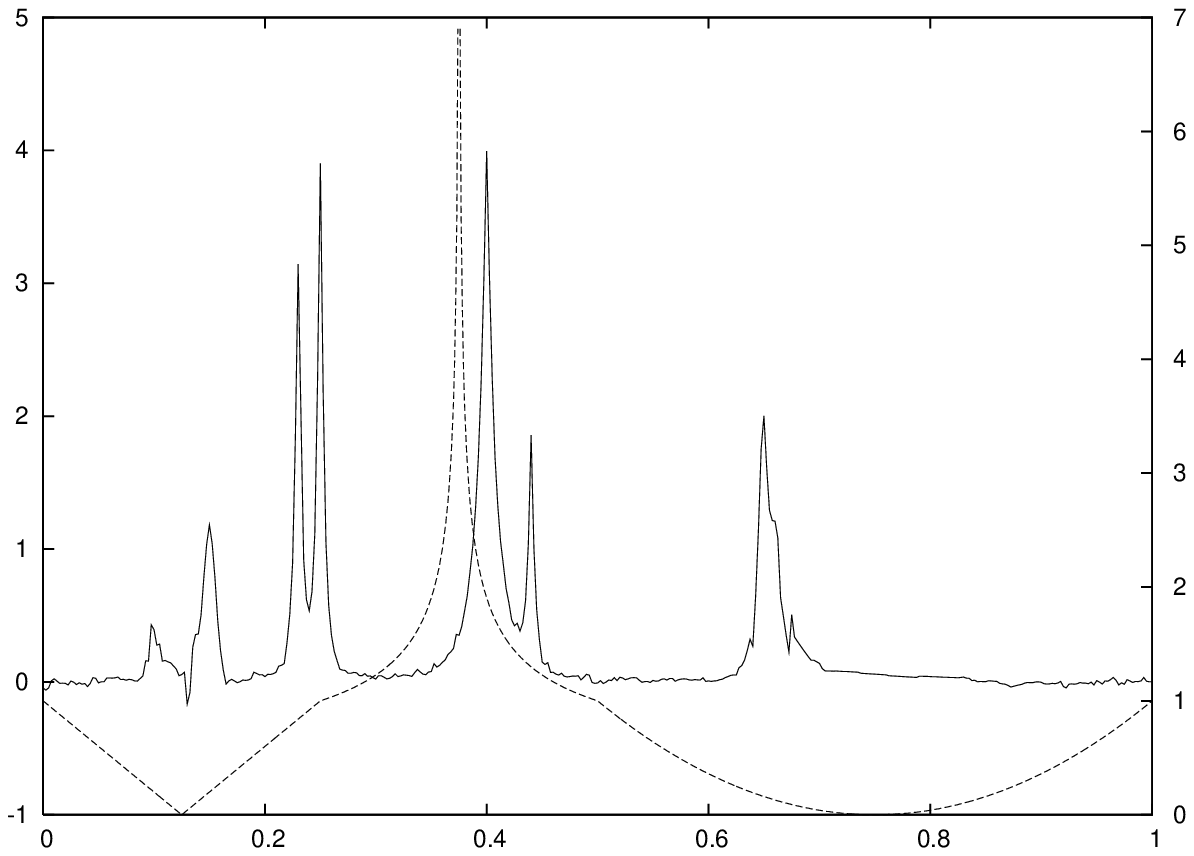}
    \includegraphics[width =
    \figurelength]{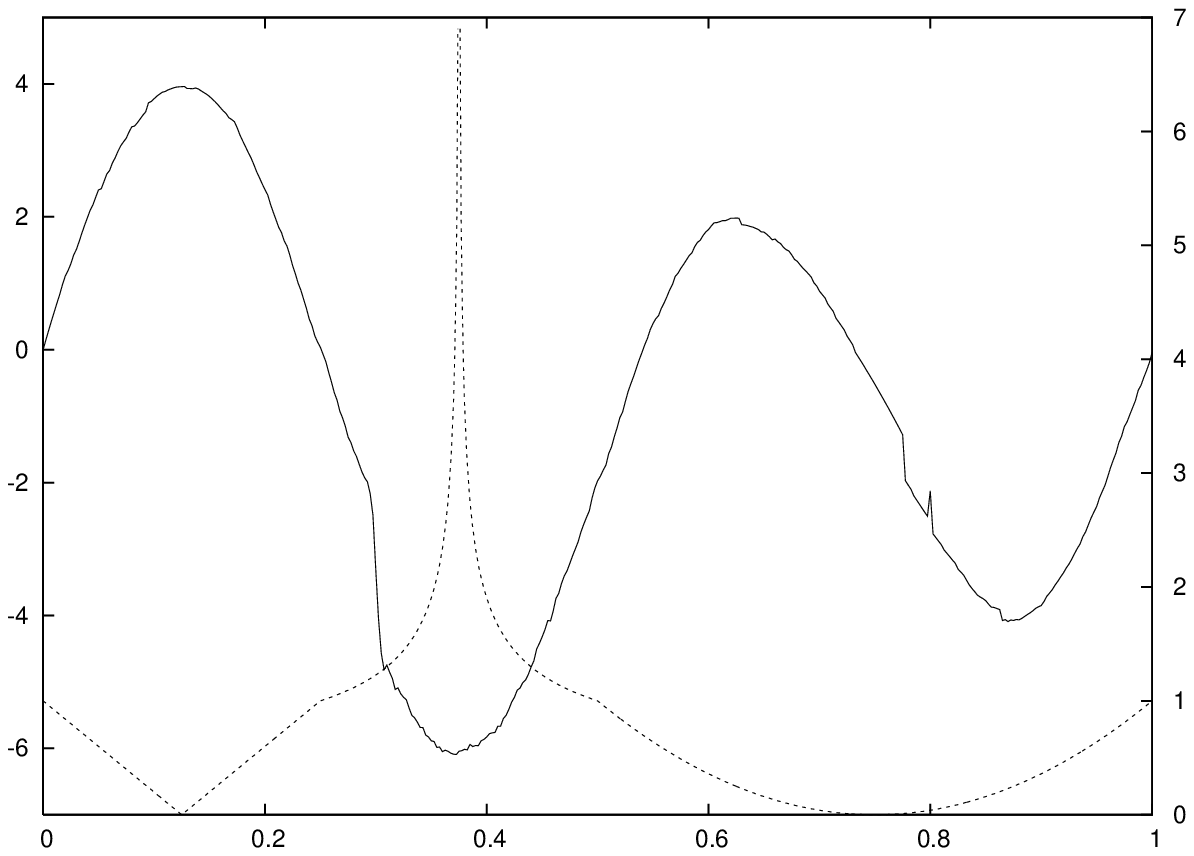}
    \includegraphics[width =
    \figurelength]{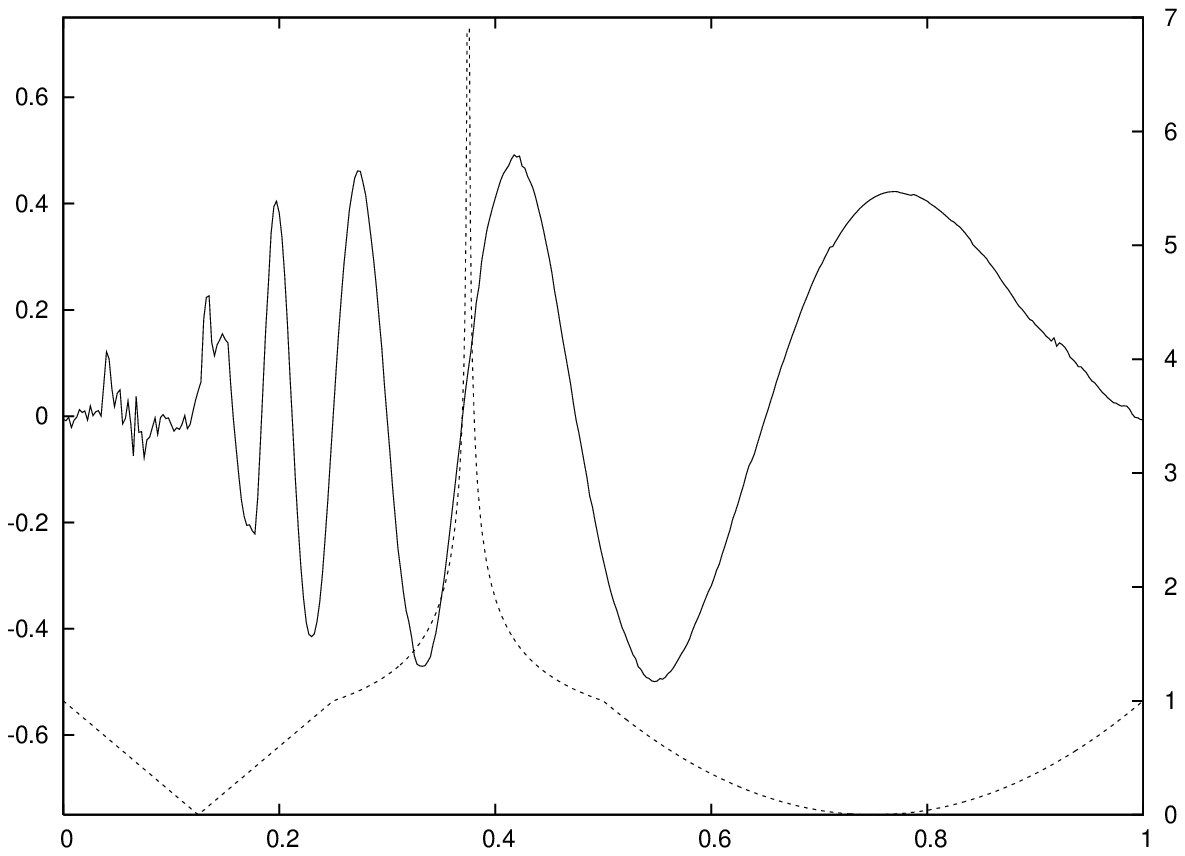}
  \end{center}
  \caption{Estimates based on the datasets in figure
    \ref{fig_datasets_example1_design}.}
  \label{fig_estimates_example1_design}
\end{figure}

Note that these estimates can be slightly improved with case by case
tuned parameters: for instance, for the first dataset (blocks), the
choice $\kpa = 0$ gives a slightly better looking estimate (the target
function is constant by parts). In figure
\ref{fig_datasets_example1_design} we show datasets with the same
signal-to-noise ratio and sample size as in figure
\ref{fig_datasets_unif_design} but the design is non-uniform (we plot
the design density on each of them). We show the estimates based on
these datasets in figure \ref{fig_estimates_example1_design}. The same
parameters as for figure \ref{fig_estimates_unif_design} are used.


\newpage

\begin{figure}[htbp]
  \begin{center}
    \includegraphics[width =
    \figurelength]{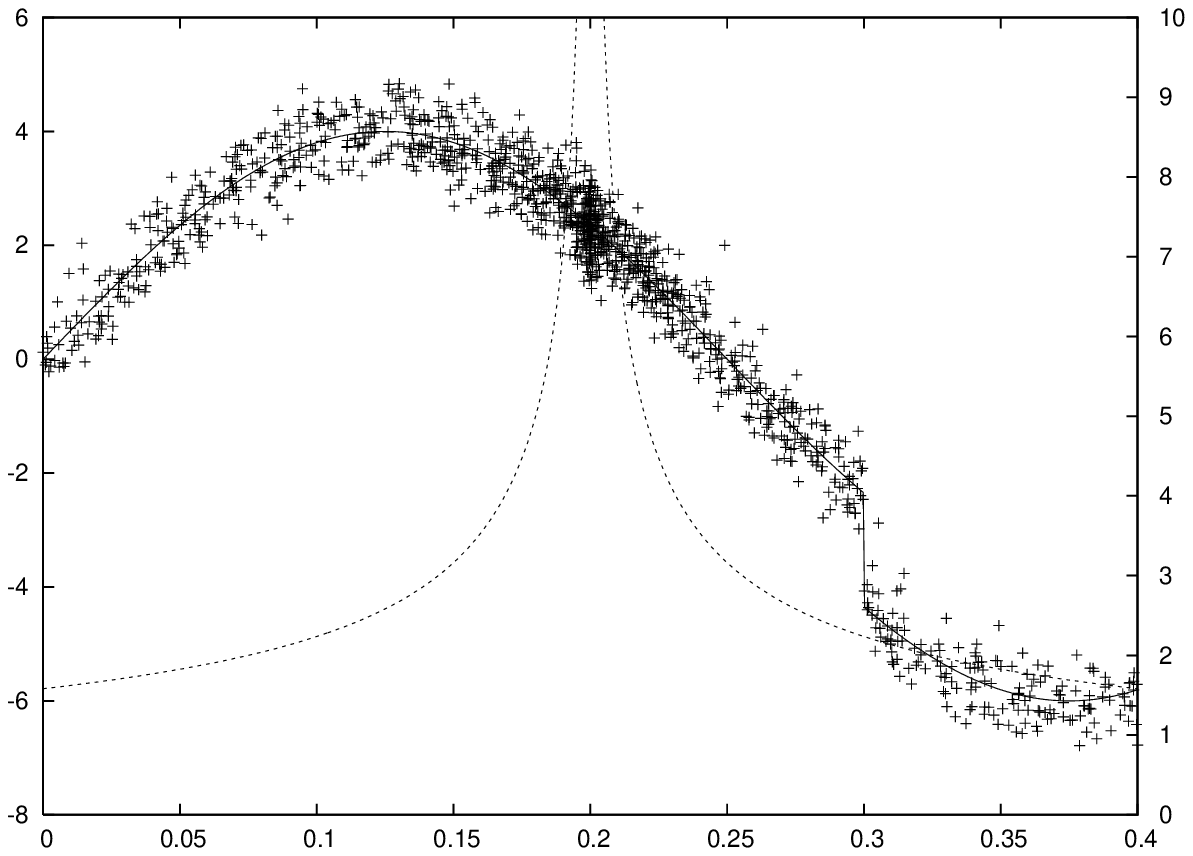}
    \includegraphics[width =
    \figurelength]{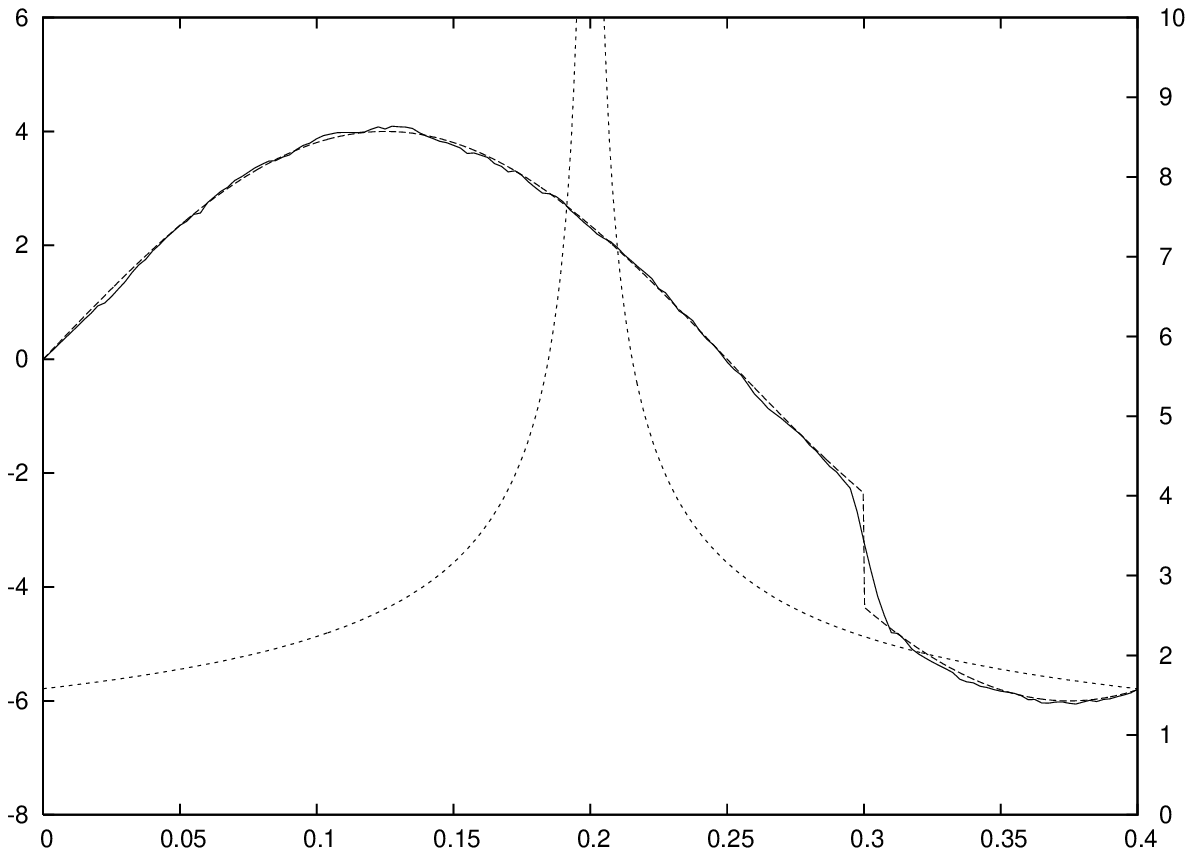}
    \includegraphics[width =
    \figurelength]{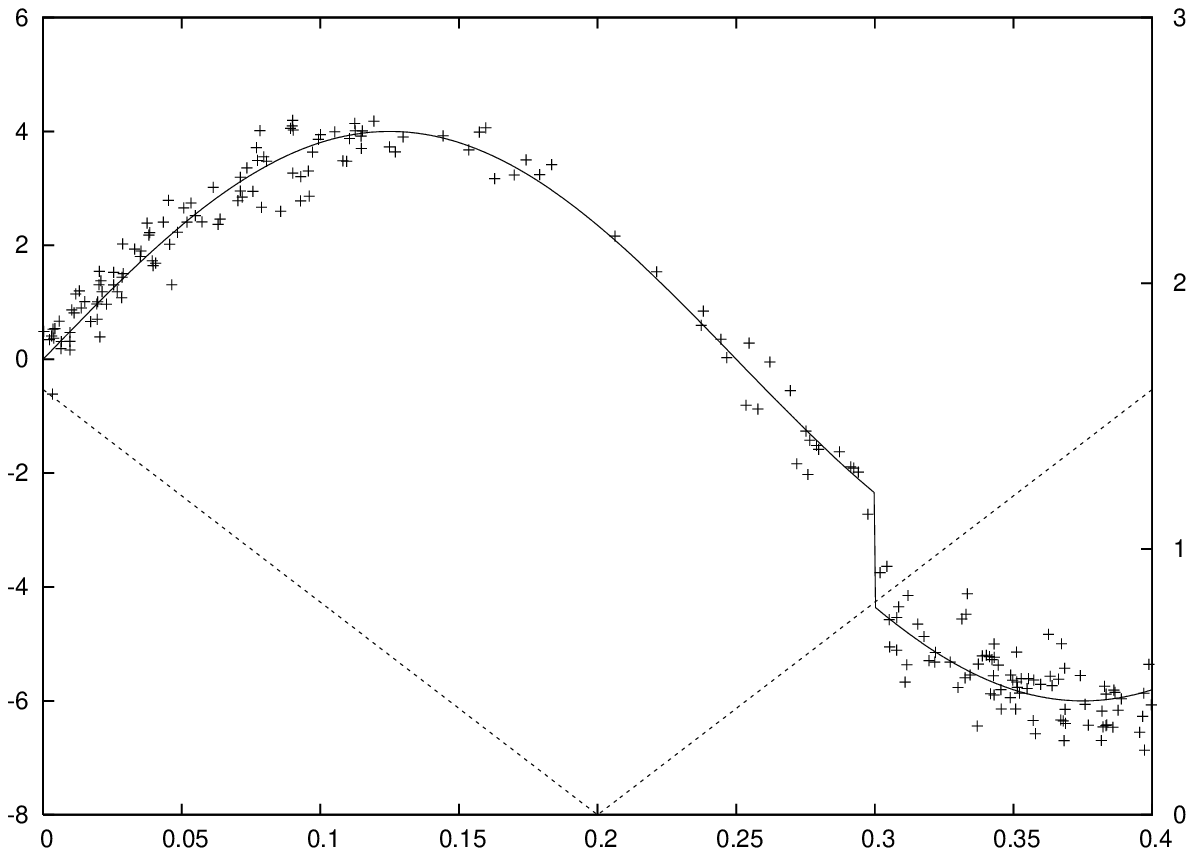}
    \includegraphics[width =
    \figurelength]{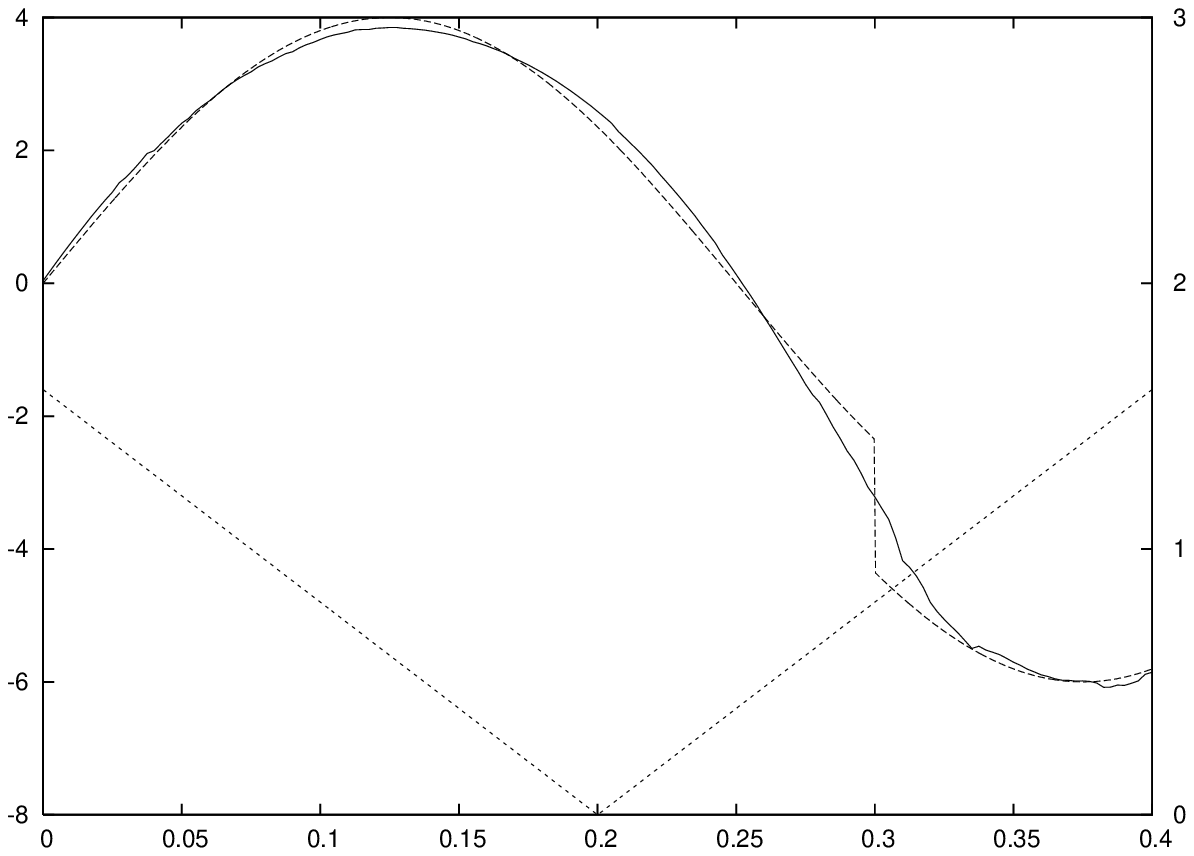}
  \end{center}
  \caption{Heavysine datasets and estimates with design density
    \eqref{eq:design_density_example} with $x_0 = 0.2$ and $\beta =
    -0.5$ at top, $\beta = 1$ at bottom.}
  \label{fig_datasets_heavysine_local1}
\end{figure}

\begin{figure}[htbp]
  \begin{center}
    \includegraphics[width =
    \figurelength]{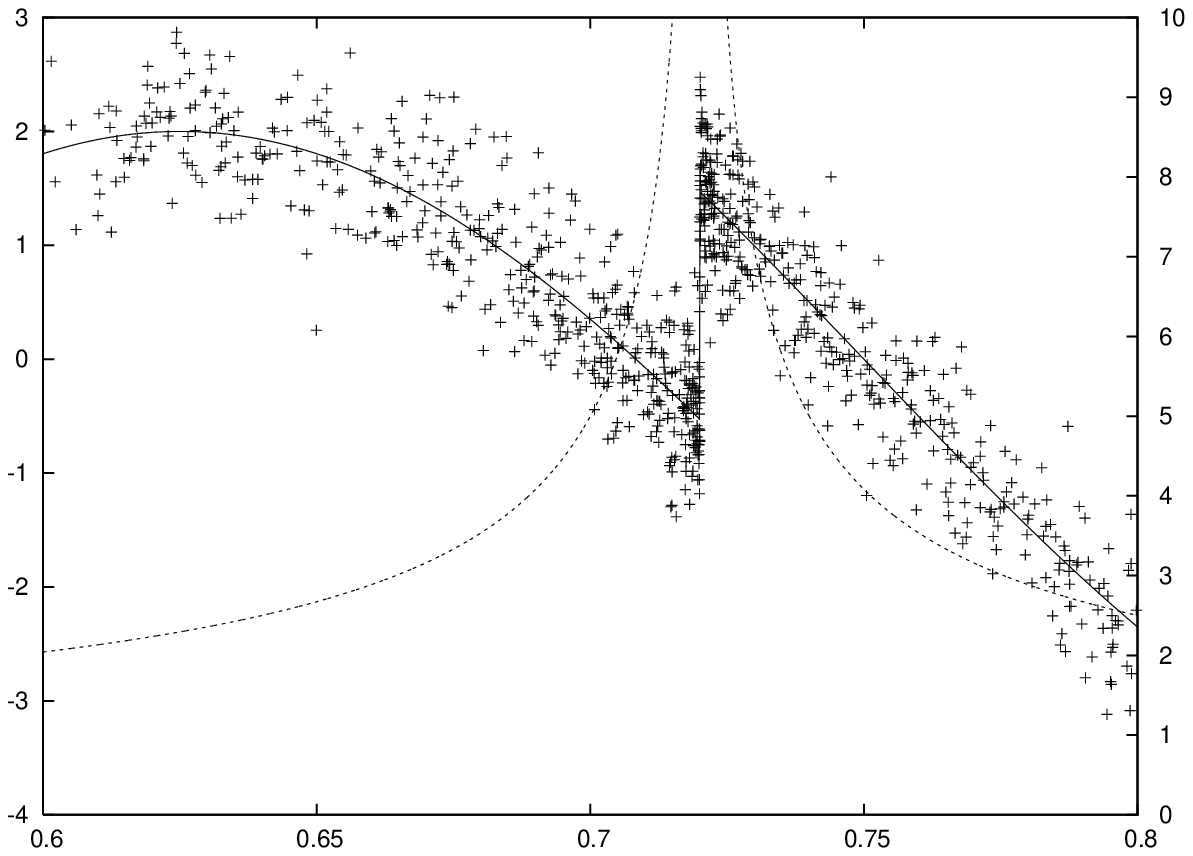}
    \includegraphics[width =
    \figurelength]{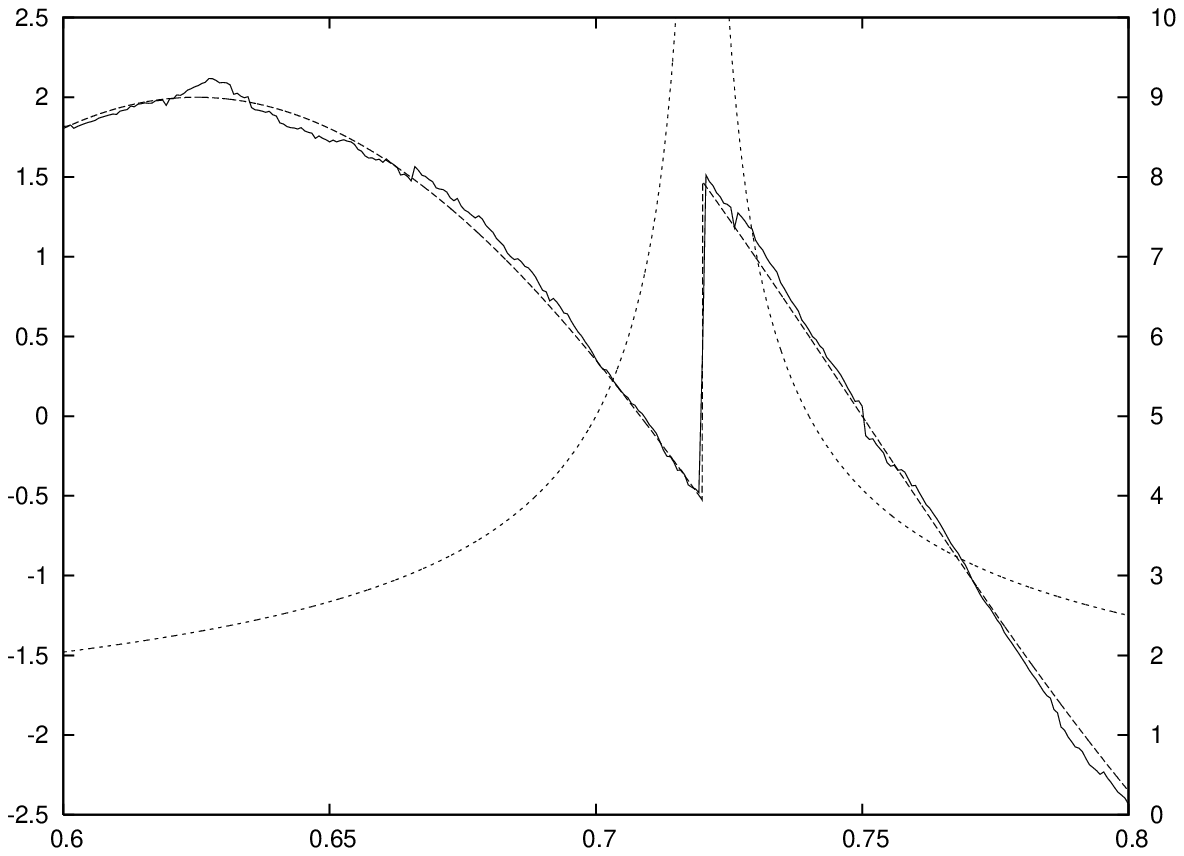}
    \includegraphics[width =
    \figurelength]{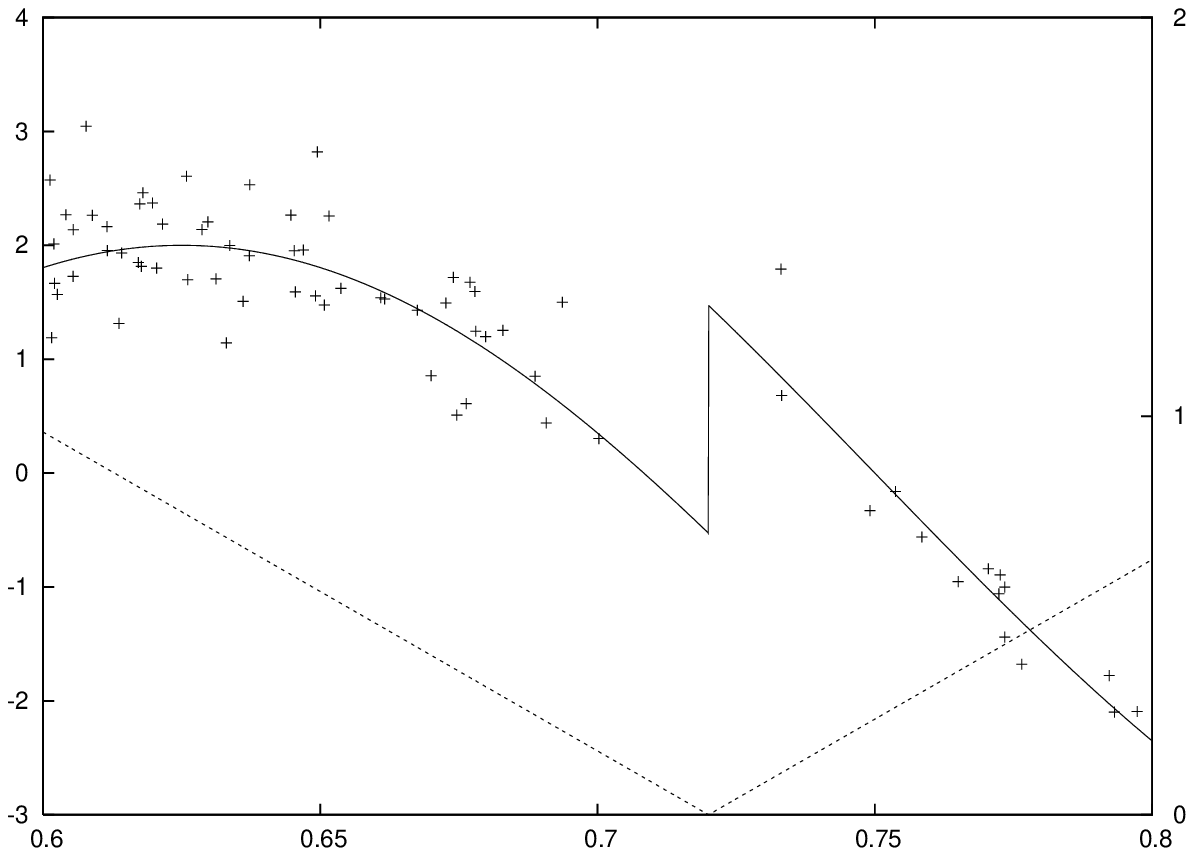}
    \includegraphics[width =
    \figurelength]{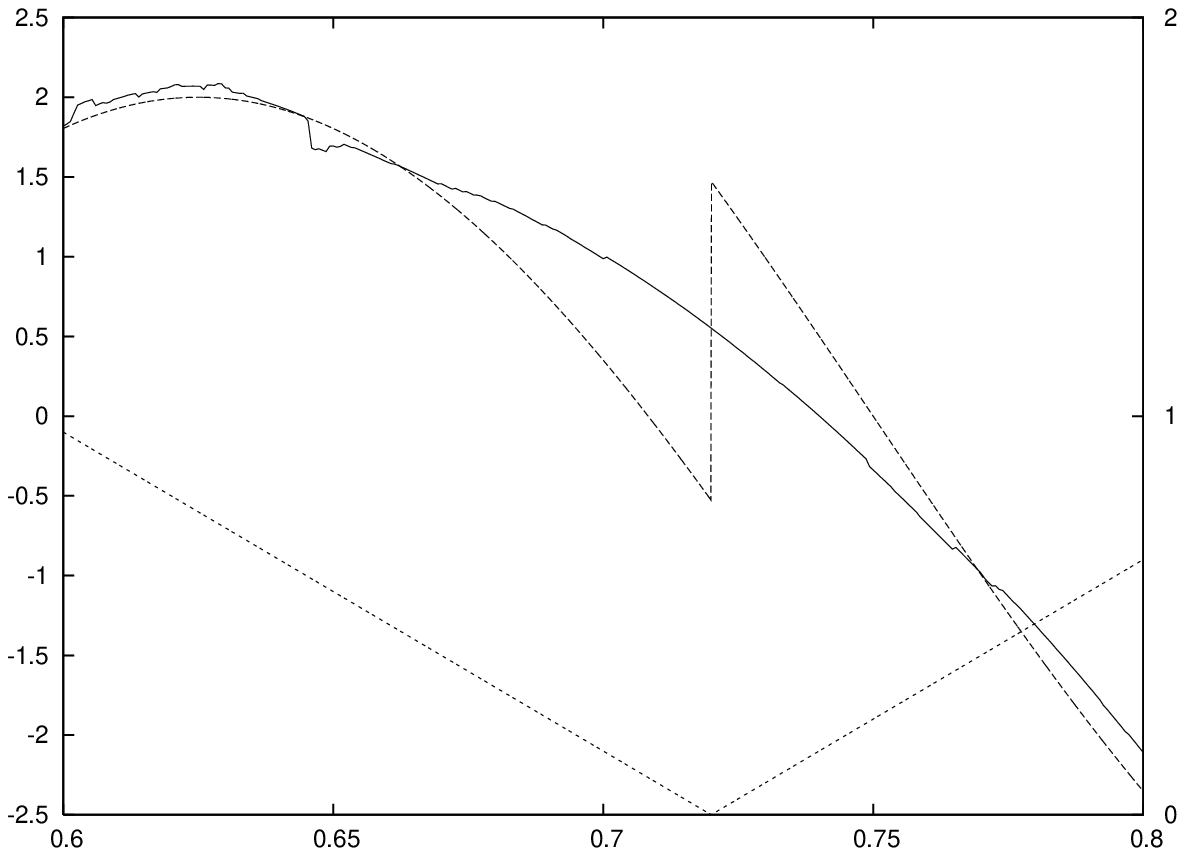}
  \end{center}
  \caption{Heavysine datasets and estimates with design density
    \eqref{eq:design_density_example} with $x_0 = 0.72$ and $\beta =
    -0.5$ at top, $\beta = 1$ at bottom.}
  \label{fig_datasets_heavysine_local2}
\end{figure}

In figures \ref{fig_datasets_heavysine_local1} and
\ref{fig_datasets_heavysine_local2} we give a more localised
illustration of the heavysine dataset. We keep the same
signal-to-noise ratio and sample size. We consider the design density
\begin{equation}
  \label{eq:design_density_example}
  \mu(x) = \frac{\beta + 1}{x_0^{\beta + 1} + (1 - x_0)^{\beta+1}}
  \bigl|x - x_0 \bigr|^{\beta} \ind{[0,1]}(x),
\end{equation}
for $x_0=0.2, 0.72$ and $\beta=-0.5, 1$.


\newpage

\section{Discussion}
\label{sec:discussion}

\subsection{On the procedure}
\label{sec:remarks_on_the_bandwidth_selection_rule}

\begin{itemize}
\item It is important to note that on the event $\Omega_h$ the
  estimator $\wh f_{h,\kpa}$ is equal to the classical local
  polynomial estimator defined by
  \begin{equation}
    \label{eq:least_square_eq_def}
    \wh f_{h,\kpa} = \arg\min_{\displaystyle g \in V_{\kpa}} \norm{g -
      Y}_h^2,
  \end{equation}
  where $V_{\kpa} = \text{Span}\{(\phi_{j})_{j=0,\ldots,\kpa}\}$. A
  necessary condition for $\wh f_{h,\kpa}$ to minimise
  \eqref{eq:least_square_eq_def} is to be solution of the linear
  problem
  \begin{equation}
    \label{eq:_def_est_variational}
    \text{find } \wh f \in V_{\kpa} \text{ such that } \forall \phi \in
    V_{\kpa}, \quad \prodscah{\wh f}{\phi} = \prodscah{Y}{\phi}. 
  \end{equation}
  The main idea of the procedure is the following: if $h$ is a
  \textit{good} bandwidth, then for any $h' \leq h$ and for all $\phi
  \in V_{\kpa}$ we should have in view of
  \eqref{eq:_def_est_variational}:
  \begin{equation*}
    \prodsca{\wh f_h - \wh f_{h'}}{\phi}_{h'} = \prodsca{\wh f_h -
      Y}{\phi}_{h'} \approx \prodsca{\xi}{\phi}_{h'},
  \end{equation*}
  which means that the difference $\wh f_h - \wh f_{h'}$ is mainly
  noise, in the sense that $\sigma^{-1} \norm{\phi}_{h'}^{-1}
  \prodsca{\wh f_h - \wh f_{h'}}{\phi}_{h'}$ is
  close in law to a standard Gaussian. \\

\item The procedure \eqref{eq:H_n_hat_def} looks like the Lepski
  procedure: in a model where the estimators can be well sorted by
  their respective variances (this is the case with kernel estimators
  in the white noise model, see Lepski and Spokoiny
  (1997)\nocite{lepski_spok97}), the Lepski procedure selects the
  largest bandwidth such that the corresponding estimator does not
  differ significantly from estimators with a smaller bandwidth.  Here
  the idea is the same, but the proposed procedure is additionally
  sensitive to the design. \\

\item The estimator $\wh f_n(x_0)$ only depends on $\kpa$ and on the
  grid $\mc H$ (to be chosen by the statistician). It does not depend
  on the regularity of $f$ nor any assumption on $\mu$. In this sense,
  this estimator is adaptive in both regularity and design. \\

\item Note that $\mb X_h = \trans \mb F_h \mb F_h$ where $\mb F_h$ is
  the matrix of size $n \times (\kpa+1)$ with entries $(\mb F_h)_{i,j}
  = (X_i-x_0)^j$ for $0 \leq i \leq n$ and $ 0 \leq j \leq \kpa$, and
  that $\ker{\mb X_h} = \ker{\mb F_h}$. Thus when $n < \kpa + 1$, $\mb
  X_h$ is not invertible since its kernel is not zero, and $\Omega_h =
  \emptyset$. This is the reason why theorem
  \ref{thm:NA_adapt_upper_bound} is stated for $n \geq \kpa+1$ and in
  the step $3$ of the procedure (see section
  \ref{sec:implementatio_of_the_procedure}) we must take $m \geq
  \kpa+1$ so that each interval in $\mc G$ contains at least $\kpa +
  1$ observations $X_i$. \\

\item The reason why we need to take the grid $\mc H = \mc
  H_1^{\tup{arith}}$ in theorem \ref{thm:asympt_adaptive_upper_bound}
  is linked with the control of $\lba_{n, \omega}$.  We can prove the
  theorem with a geometrical grid if we additionally assume $\lba_{n,
    \omega} > \lba$ for $\lba > 0$, but we preferred to work only
  under the regularly varying design assumption with a restricted grid
  choice without extra assumption on the model. \\

\item The fact that the noise level $\sigma$ is known is of little
  importance. If it is unknown we can plug-in some estimator $\wh
  \sigma_n^2$ in place of $\sigma^2$. Following Gasser \etal
  (1986)\nocite{gasser_sroka86} or Buckley \etal
  (1988)\nocite{buckley_eagleson88} we can consider
  \begin{equation}
    \label{eq:sigma_estim}
    \wh \sigma_n^2 = \frac{1}{2(n-1)} \sum_{i=1}^{n-1}(Y_{(i+1)} -
    Y_{(i)})^2,
  \end{equation}
  where $Y_{(i)}$ is the observation at the point $X_{(i)}$ where
  $X_{(1)} \leq X_{(2)} \leq \ldots \leq X_{(n)}$. 
\end{itemize}

\subsection{Comparison with previous results}
\label{sec:comparison_with_existing_results}

\begin{itemize}
\item In Guerre (1999)\nocite{guerre99}, for the estimation of the
  regression function at the point $0$ in a more general setup for the
  design, the author works conditionally on $\mf X_n$ and gives an
  upper bound with a data-driven rate similar to
  \eqref{eq:R_n_omega_def}. The author considers then as an example
  the case of an i.i.d. design with density $\mu$ such that $\mu(x)
  \sim x^{\beta}$ close to $0$ for $\beta > -1$, which is a particular
  case of regularly varying density at $0$ of index $\beta$. Here the
  approach is the same: under the regular variation assumption we
  derive from theorem \ref{thm:NA_adapt_upper_bound} an asymptotic
  upper-bound with a deterministic rate (theorem
  \ref{thm:asympt_adaptive_upper_bound}). \\

\item Bandwidth selection procedures in local polynomial estimation
  can be found in Fan and Gijbels (1995)\nocite{fan_gijbels95},
  Goldenshluger and Nemirovski
  (1997)\nocite{goldenshluger_nemirovski97} or Spokoiny
  (1998)\nocite{spok98}. In this last paper the author is interested
  in the regression function estimation near a change point. The main
  idea and difference between the work by Spokoiny
  (1998)\nocite{spok98} and the previous work by Goldenshluger and
  Nemirovski (1997)\nocite{goldenshluger_nemirovski97} is to solve the
  linear problem \eqref{eq:_def_est_variational} in a non symmetrical
  neighbourhood of $x_0$ not containing the change point.  Our
  adaptive procedure \eqref{eq:H_n_hat_def} is mainly inspired from
  the work of Spokoiny and adapted for the degenerate random design
  problem. We have also made improvements, for instance we do not need
  to bound the estimator and the function at $x_0$ by some known
  constant.
\end{itemize}

\section{Proofs}
\label{sec:proofs}

In the following we denote by $\mb P_{k,h}$ the projection in the
space $V_k$ for the scalar product $\prodscah{\cdot}{\cdot}$. We
denote respectively by $\prodsca{\cdot}{\cdot}$ and by $\norm{\cdot}$
the Euclidean scalar product and the Euclidean norm in
$\setR^{\kpa+1}$.  We denote by $\norminfty{\cdot}$ the sup norm in
$\setR^{\kpa+1}$. We define $e_1 \eqdef (1,0,\ldots,0)$, the first
canonical basis vector in $\setR^{\kpa+1}$. 

\subsection{Preparatory results and proof of theorem
  \ref{thm:NA_adapt_upper_bound}}

The next lemma is a version of the local polynomial estimator
bias-variance decomposition, which is classical: see Cleveland
(1979)\nocite{cleveland79}, Tsybakov (1986)\nocite{tsybakov86},
Korostelev and Tsybakov (1993)\nocite{korostelev_tsybakov93}, Fan and
Gijbels (1995, 1996)\nocite{fan_gijbels95,fan_gijbels96},
Goldenshluger and Nemirovski
(1997)\nocite{goldenshluger_nemirovski97}, Spokoiny
(1998)\nocite{spok98} and Tsybakov (2003)\nocite{tsybakov03}, among
others. The version given by lemma \ref{lem:bias_variance_adapt} is
close to the one in Spokoiny (1998)\nocite{spok98}.  Let us introduce
for any positive integer $k$ the continuity modulus
\begin{equation*}
  \omega_{f,k}(x_0,h) = \inf_{P \in \mc P_{k}} \sup_{|x-x_0| \leq h}
  |f(x) - P(x-x_0)|. 
\end{equation*}
Note that if $k_1 \leq k_2$ we clearly have $\omega_{f, k_2}(x_0,h)
\leq \omega_{f,k_1}(x_0,h)$.

\begin{lemma}[Bias variance decomposition]
  \label{lem:bias_variance_adapt}
  On the event $\Omega_h$ the estimator $\wh f_{h,\kpa}$ from
  definition \ref{def:loc_pol_est_final} satisfies for any $k \leq
  \kpa$
  \begin{equation}
    \label{eq:bs_pl2}
    |\wh f_h(x_0) - f(x_0)| \leq \lambda^{-1}(\mc G_h) \sqrt{\kpa+1}
    \bigl( \omega_{f,k}(x_0,h) + \sigma N_{n,h}^{-1/2} |\gamma_h|
    \bigr),
  \end{equation}
  where $\gamma_h$ is, conditional on $\mf X_n$, centered Gaussian
  such that $\Efm\{\gamma_h^2 | \mf X_n\} \leq 1$. 
\end{lemma}

\begin{proof}
  On $\Omega_h$ we have $\mbt X_h = \mb X_h$ and $\lba(\mb X_h) >
  N_{n,h}^{-1/2} > 0$, then $\mb X_h$ is invertible. Since $\Lba_h$ is
  clearly invertible on this event, $\mc G_h$ is also invertible. Let
  $0 < \von \leq \frac{1}{2}$. By definition of $\omega_{f, \kpa}(x_0,
  h)$ we can find a polynomial $P_{f,h}^{\von} \in \mc P_{\kpa}$ such
  that
  \begin{equation*}
    \sup_{x \in [x_0 - h, x_0 + h]} |f(x) - P_{f,h}^{\von}(x)| \leq
    \omega_{f,\kpa}(x_0,h) + \frac{\von}{\sqrt{n}}. 
  \end{equation*}
  In particular we have $| f(x_0) - P_{f,h}^{\von}(x_0)| \leq
  \frac{\von}{\sqrt{n}}$ and if we denote by $\tta_{h}$ the
  coefficients vector of $P_{f,h}^{\von}$ then
  \begin{equation*}
    |\wh f_{h,\kpa}(x_0) - f(x_0)| \leq |\prodsca{\Lba_h^{-1} (\wh \tta_h -
      \tta_{h})}{e_1}| + \frac{\von}{\sqrt{n}} = |\prodsca{\mc G_h^{-1}
      \Lba_h \mb X_h (\wh \tta_h - \tta_{h})}{e_1}| +
    \frac{\von}{\sqrt{n}} . 
  \end{equation*}
  Then in view of \eqref{eq:_def_est_variational} one has for $j =
  0,\ldots,\kpa$:
  \begin{align*}
    (\mb X_h (\wh \tta_h - \tta_{h}))_j = \prodscah{\wh f_{h,\kpa} -
      P_{f, h}^{\von}}{\phi_j} = \prodscah{Y - P_{f,h}^{\von}}{\phi_j}
    = \prodscah{f - P_{f,h}^{\von}}{\phi_j} + \prodscah{\xi}{\phi_j},
  \end{align*}
  thus we can decompose $\mb X_h (\wh \tta_h - \tta_{h}) \eqdef B_h +
  V_h$ and then:
  \begin{equation*}
    |\wh f_{h,\kpa}(x_0) - f(x_0)| \leq |\prodsca{\mc G_h^{-1} \Lba_h
      B_h}{e_1}| + |\prodsca{\mc G_h^{-1} \Lba_h V_h}{e_1}| +
    \frac{\von}{\sqrt{n}} \eqdef A + B + \frac{\von}{\sqrt{n}}. 
  \end{equation*}
  We have
  \begin{align*}
    A &\leq \norm{\mc G_h^{-1} \Lba_h B_h} \leq \norm{\mc G_h^{-1}}
    \norm{\Lba_h B_h} \leq \norm{\mc G_h^{-1}} \sqrt{\kpa+1}
    \norminfty{\Lba_h B_h},
  \end{align*}
  and
  \begin{equation*}
    |(\Lba_h B_h)_j| = \normh{\phi_j}^{-1} |\prodscah{f -
      P_{f,h}^{\von}}{\phi_j} | \leq \normh{f - P_{f,h}^{\von}}
    \leq \omega_{f,\kpa}(x_0,h) + \frac{\von}{\sqrt{n}}. 
  \end{equation*}
  For any symmetrical and positive matrix $M$ we have $\lba^{-1}(M) =
  \norm{M^{-1}}$ then since $\norm{\Lba_h^{-1}} \leq 1$ we have on the
  event $\Omega_h$:
  \begin{equation*}
    \norm{\mc G_h^{-1}} = \norm{\Lba_h^{-1} \mb X_h^{-1} \Lba_h^{-1}}
    \leq \norm{\mb X_h^{-1}} = \lba^{-1}(\mb X_h) \leq N_{n,h}^{1/2}
    \leq \sqrt{n}. 
  \end{equation*}
  Thus $A \leq \norm{\mc G_h^{-1}} \sqrt{\kpa+1} \omega_{f, \kpa}(x_0,
  h) + \von \sqrt{\kpa+1} \leq \norm{\mc G_h^{-1}} \sqrt{\kpa+1}
  \omega_{f, k}(x_0, h) + \von \sqrt{\kpa+1}$ since $k \leq \kpa$. 
  Conditional on $\mf X_n$, the random vector $V_h$ is centered
  Gaussian with covariance matrix $\sigma^2 N_{n,h}^{-1} \mb X_h$. 
  Thus $\mc G_h^{-1} \Lba_h V_h$ is again centered Gaussian, with
  covariance matrix
  \begin{equation*}
    \sigma^2 N_{n,h}^{-1} \mc G_h^{-1} \Lba_h \mb X_h \Lba_h \mc
    G_h^{-1} = \sigma^2 N_{n,h}^{-1} \mc G_h^{-1},
  \end{equation*}
  and $B$ is then centered Gaussian with variance
  \begin{equation*}
    \sigma^2 N_{n,h}^{-1} \prodsca{e_1}{\mc G_h^{-1} e_1} \leq \sigma^2
    N_{n,h}^{-1} \norm{\mc G_h^{-1}}. 
  \end{equation*}
  Since $\mc G_h$ is positive symmetrical and its entries are smaller
  than one in absolute value we get $\norm{\mc G^{-1}} = \lba^{-1}(\mc
  G_h)$ and $\lba(\mc G_h) = \inf_{\norm{x} = 1} \prodsca{x}{\mc G_h
    x} \leq \norm{\mc G_h e_1} \leq \sqrt{\kpa+1}$.  Thus $\norm{\mc
    G_h^{-1}} \leq \sqrt{\kpa+1} \norm{\mc G_h^{-1}}^2$, and the
  proposition follows. 
\end{proof}

Let us introduce the events 
\begin{equation*}
  \mc A_{h',h, j} \eqdef \bigl\{|\prodscahp{\wh f_{h, \kpa} - \wh
    f_{h',\kpa}}{\phi_j}| \leq \sigma \norm{\phi_j}_{h'} T_{n,h',h}
  \bigr\},
\end{equation*}
$\mc A_{h',h} \eqdef \bigcap_{j=0}^{\kpa} \mc A_{h', h, j}$ and $ \mc
A_h \eqdef \bigcap_{h' \in \mc H_h} \mc A_{h',h}$. The following lemma
shows that if some bandwidth $h$ is \emph{good} in the sense that $h
\leq H_{n,\omega}$ ($h$ is smaller than the ideal adaptive bandwidth)
then $h$ can be selected by the procedure with a large probability.

\begin{lemma}
  \label{lem:prob_bandwidth_reject}
  Let $f \in \mc F_{\delta}(x_0, \omega)$ for $\omega \in \RV(s)$ with
  $0 < s \leq \kpa + 1$. If $h$ is such that $h \leq H_{n,\omega}
  \wedge \delta$ we have on $\Omega_h$ for any $n \geq \kpa+1$\tup:
  \begin{equation*}
    \Pfm \bigl\{ \mc A_{h} | \mf X_n \bigr\} \geq 1 - (\kpa+1)
    N_{n,h}^{-2 p}. 
  \end{equation*}
\end{lemma}

\begin{proof}
  Let $j \in \{ 0,\ldots, \kpa\}$ and $h' \in \mc H_h$. On $\Omega_h$
  we have in view of \eqref{eq:least_square_eq_def} that $\wh
  f_{h,\kpa} = \mb P_{\kpa,h}(Y)$ thus using
  \eqref{eq:_def_est_variational} we can decompose:
  \begin{align*}
    \prodsca{\wh f_{h', \kpa} - \wh f_{h, \kpa}}{\phi_j}_{h'} =
    \prodsca{Y - \wh f_{h,\kpa}}{\phi_j}_{h'} &= \prodsca{f - \wh
      f_{h, \kpa}}{\phi_j}_{h'} + \prodsca{\xi}{\phi_j}_{h'} \\
      &= \prodsca{f - \mb P_{\kpa,h}(f)}{\phi_j}_{h'} + \prodsca{\mb
        P_{\kpa,h}(f) - \wh f_{h, \kpa} }{\phi_j}_{h'} +
      \prodsca{\xi}{\phi_j}_{h'} \\
      &= \prodsca{f - \mb P_{\kpa,h}(f)}{\phi_j}_{h'} + \prodsca{\mb
        P_{\kpa,h}(f - Y)}{\phi_j}_{h'} + \prodsca{\xi}{\phi_j}_{h'} \\
      &= \prodsca{f - \mb P_{\kpa,h}(f)}{\phi_j}_{h'} - \prodsca{\mb
        P_{\kpa,h}(\xi)}{\phi_j}_{h'} + \prodsca{\xi}{\phi_j}_{h'} \\
      &\eqdef A + B + C. 
  \end{align*}
  The term $A$ is a bias term. By the definition of
  $\omega_{f,k}(x_0,h)$ we can find a polynomial $P_{f, h}^n \in
  V_{k}$ such that
  \begin{equation*}
    \sup_{x \in [x_0 - h, x_0 + h]} |f(x) - P_{f, h}^n(x)| \leq
    \omega_{f,k}(x_0,h) + \von_n,
  \end{equation*}
  where $\von_n \eqdef \frac{C_{\kpa} \sigma}{2} \sqrt{\frac{C_{p}
      \log 2}{n}}$ (see \eqref{eq:T_n_def}). Since $h' \leq h \leq
  \delta$, $f \in \mc F_{\delta}(x_0,\omega)$ and $P_{f, h}^n \in V_k
  \subset V_{\kpa}$ we get
  \begin{align*}
    |A| \leq \norm{f - \mb P_{\kpa,h}(f)}_{h'} \norm{\phi_j}_{h'}
    &\leq \normh{f - P_{f, h}^n - \mb P_{\kpa, h}(f - P_{f, h}^n)}
    \norm{\phi_j}_{h'} \\
    &\leq \normh{f - P_{f, h}^n} \norm{\phi_j}_{h'} \\
    &\leq \norm{\phi_j}_{h'} (\omega_{f,k}(x_0,h) + \von_n) \leq
    \norm{\phi_j}_{h'} ( \omega(h) + \von_n),
  \end{align*}
  since $\mb P_{\kpa,h}$ is a projection with respect to
  $\prodscah{\cdot}{\cdot}$. If $h < H_{n, \omega}$ we have in view of
  \eqref{eq:H_n_def_adapt} that $\omega(h) \leq \sigma \sqrt{N_{n,
      h}^{-1} \log n}$. When $h = H_{n, \omega}$ two cases can occur. 
  If the graphs of $h \mapsto \sigma \sqrt{N_{n,h}^{-1} \log n}$ and
  $h \mapsto \omega(h)$ cross each other we have $\omega(h) = \sigma
  \sqrt{N_{n,h}^{-1} \log n}$. When these graphs do not cross we
  introduce $H_{n, \omega}^- = \max\{ h \in \mc H | h < H_{n, \omega}
  \}$ and $H_{n, \omega}^+ = \min\{ h \in \mc H | h \geq H_{n, \omega}
  \}$.  If $\mc H = \agrid$ we have $N_{n, H_{n, \omega}} \leq N_{n,
    H_{n, \omega}^+} \leq N_{n,H_{n, \omega}^-} + a$ while when $\mc H
  = \ggrid$ we get $N_{n,H_{n, \omega}} \leq N_{n, H_{n, \omega}^+}
  \leq (1+a) N_{n, H_{n, \omega}^-}$.  Then for any $h \leq H_{n,
    \omega}$:
  \begin{equation}
    \label{eq:A_estimate}
    |A| \leq
    \begin{cases}
      \norm{\phi_j}_{h'} (\sigma \sqrt{(N_{n,h} - a)^{-1} \log n} +
      \von_n) &\text{ if } \mc H = \agrid, \\
      \norm{\phi_j}_{h'} (\sigma \sqrt{(1 + a) N_{n,h}^{-1} \log n} +
      \von_n) &\text{ if } \mc H = \ggrid. 
    \end{cases}
  \end{equation}
  Conditional on $\mf X_n$, $B$ and $C$ are centered Gaussian. We have
  $\mc L(C | \mf X_n) = \mc N(0, \sigma^2 N_{n,h'}^{-1}
  \norm{\phi_j}_{h'}^2)$ and conditional on $\mf X_n$ the vector $\mb
  P_{\kpa,h}(\xi)$ is centered Gaussian with covariance matrix
  \begin{equation*}
    \sigma^2 \mb P_{\kpa,h} \trans \mb P_{\kpa,h} = \sigma^2
    \mb P_{\kpa,h},
  \end{equation*}
  since $\mb P_{\kpa,h}$ is a projection. Thus $B$ is centered
  Gaussian with variance
  \begin{align*}
    \Efm\{ \prodsca{\mb P_{\kpa,h}(\xi)}{\phi_j}_{h'}^2 | \mf X_n \}
    &\leq \norm{\phi_j}_{h'}^2 \Efm\{ \norm{\mb
      P_{\kpa,h}(\xi)}_{h'}^2 | \mf X_n \}  \\
    &= N_{n,h'}^{-1} \norm{\phi_j}_{h'}^2 \text{tr}(\Var(\mb
    P_{\kpa,h}(\xi)|\mf X_n )) \\
    &= \sigma^2 N_{n,h'}^{-1} \norm{\phi_j}_{h'}^2 \text{tr}(\mb
    P_{\kpa,h}) \\
    &\leq \sigma^2 N_{n,h'}^{-1} \norm{\phi_j}_{h'}^2 \dim(V_{\kpa})
    \leq \sigma^2 N_{n,h'}^{-1} \norm{\phi_j}_{h'}^2 (\kpa+1),
  \end{align*}
  where we last used that $\mb P_{\kpa,h}$ is the projection in
  $V_{\kpa}$. Then conditional on $\mf X_n$, $B + C$ is centered
  Gaussian with variance
  \begin{align*}
    \Efm \{ (B+C)^2 | \mf X_n \} &\leq \Efm\{B^2 + 2BC + C^2|\mf X_n
    \} \\
    &\leq \Efm\{B^2|\mf X_n \} + 2\sqrt{\Efm\{B^2|\mf X_n \}\Efm\{
      C^2|\mf X_n \}} + \Efm\{C^2|\mf X_n \} \\
    &\leq \sigma^2 (1 + \sqrt{\kpa + 1})^2 N_{n, h'}^{-1}
    \norm{\phi_j}_{h'}^2 C_{\kpa}^2. 
  \end{align*}
  Using \eqref{eq:A_estimate} and since $2 \leq N_{n,h} \leq n$ on
  $\Omega_h$ we have
  \begin{equation*}
    \mc A_{h',h,j}^c \subset \Bigl\{ \frac{|B + C|}{\sigma N_{n,h'}^{-1/2}
      \norm{\phi_j}_{h'} C_{\kpa} } > \sqrt{ C_{p} \log
      N_{n,h}} / 2 \Bigr\},
  \end{equation*}
  and using a standard Gaussian large deviation inequality we get
  \begin{align*}
    \Pfm\bigl\{ \mc A_{h',h,j}^c \bigl| \mf X_n \bigr\} \leq \exp
    \bigl( -(1 + 2 p) \log N_{n,h} \bigr) = N_{n,h}^{-(1+ 2p)}. 
  \end{align*}
  Since $\#(\mc H_h) \leq N_{n,h}$ we finally have
  \begin{equation*}
    \Pfm\{ \mc A_h^c | \mf X_n \} \leq \Pfm\Bigl\{ \bigcup_{h' \in \mc
      H_h} \bigcup_{j=0}^{\kpa} \mc A_{h',h,j}^c \bigl| \mf X_n
    \Bigr\} \leq (\kpa+1) N_{n,h}^{-2 p}. \qedhere
  \end{equation*}
\end{proof}

\begin{lemma}
  \label{lem:estimator_difference}
  Let $h \in \mc H$ and $h' \in \mc H_h$. On the event $\Omega_{h'}
  \cap \mc A_{h',h}$ one has\tup:
  \begin{equation*}
    |\wh f_h(x_0) - \wh f_{h'}(x_0)| \leq C_{p,\kpa,a} \norm{\mc
      G_{h'}^{-1}} \sigma \sqrt{N_{n, h'}^{-1} \log n},
  \end{equation*}
  where $C_{p, \kpa, a} \eqdef \sqrt{\kpa+1} (\sqrt{1+a} + C_{\kpa}
  \sqrt{C_{p}})$. 
\end{lemma}

\begin{proof}
  In view of definition \ref{def:loc_pol_est_final} and since $\mc
  G_{h'}$ is invertible on $\Omega_{h'}$ we have
  \begin{align*}
    |\wh f_h(x_0) - \wh f_{h'}(x_0)| = |\prodsca{\Lba_{h'}^{-1}(\wh
      \tta_h - \wh \tta_{h'})}{e_1}| &\leq \norm{\Lba_{h'}^{-1}(\wh
      \tta_h - \wh \tta_{h'})} \\
    &= \norm{\mc G_{h'}^{-1} \Lba_{h'} \mb X_{h'} (\wh \tta_h - \wh
      \tta_{h'})} \\
    &\eqdef \norm{\mc G_{h'}^{-1} \Lba_{h'} D_{h',h} } \leq \norm{\mc
      G_{h'}^{-1}} \sqrt{\kpa+1} \norminfty{\Lba_{h'} D_{h',h}}. 
  \end{align*}
  On $\mc A_{h',h}$ we have for any $j \in \{0, \ldots, \kpa\}$:
  \begin{equation*}
    |(D_{h',h})_j| = |\prodsca{\wh f_h - \wh f_{h'}}{\phi_j}_{h'}|
    \leq \sigma \norm{\phi_j}_{h'} T_{n,h',h},
  \end{equation*}
  thus $\norminfty{\Lba_{h'} D_{h',h}} \leq \sigma T_{n,h',h}$. Since
  $h' \leq h$ and $N_{n,h} \leq n$ we have when $\mc H = \ggrid$
  \begin{equation}
    \label{eq:T_n_estimate_proof_prop}
    T_{n,h',h} \leq (C_{\kpa} \sqrt{C_{p}} + \sqrt{1 + a}) \sqrt{N_{n,
        h'} \log n},
  \end{equation}
  and when $\mc H = \agrid$ we have by construction $N_{n,h} \geq 1 +
  a$ thus $(N_{n,h} - a)^{-1} \leq (1 + a) N_{n,h}^{-1}$ and
  \eqref{eq:T_n_estimate_proof_prop} holds again. 
\end{proof}

\begin{lemma}
  \label{lem:bad_case}
  For any $p, \alpha > 0$ and $0 < h' \leq h \leq 1$ the estimator
  $\wh f_{h'}$ given by definition \ref{def:loc_pol_est_final}
  satisfies\tup:
  \begin{equation*}
    \sup_{\displaystyle f \in \mc U(\alpha) } \Efm \bigl\{
    |\wh f_{h'}(x_0)|^p | \mf X_n \bigr\} \leq 
    C_{\sigma, p, \kpa} (\alpha \vee 1)^p N_{n,h}^{p / 2},
  \end{equation*}
  where $C_{\sigma, p, \kpa} = (\kpa + 1)^{p / 2} \sqrt{ \frac{2}{\pi}
  } \int_{\setR^+} (1 + \sigma t)^p \exp(-t^2 / 2) dt$. 
\end{lemma}

\begin{proof}
  If $N_{n,h'} = 0$ we have $\wh f_{h'} = 0$ and the result is
  obvious, thus we assume $N_{n,h'} > 0$. Since $\lba(\mbt X_{h'})
  \geq N_{n,h'}^{-1/2} > 0$, $\mbt X_{h'}$ and $\Lba_{h'}$ are
  invertible and also $\mc G_{h'}$. Thus,
  \begin{equation*}
    \wh f_{h'}(x_0) = \prodsca{\Lba_{h'}^{-1} \wh \tta_{h'}}{e_1} =
    \prodsca{\mc G_{h'}^{-1} \Lba_{h'} \mbt X_{h'} \wh \tta_{h'}}{e_1} =
    \prodsca{\mc G_{h'}^{-1} \Lba_{h'} \mb Y_{h'}}{e_1}. 
  \end{equation*}
  For any $j \in \{ 0,\ldots,\kpa \}$ we have $(\Lba_{h'} \mb
  Y_{h'})_j = \normhp{\phi_j}^{-1} ( \prodscahp{f}{\phi_j} +
  \prodscahp{\xi}{\phi_j}) \eqdef B_{h',j} + V_{h',j}$. Since $f \in
  \mc U(\alpha)$ we have
  \begin{equation*}
    |B_{h',j}| \leq \normhp{\phi_j}^{-1} |\prodscahp{f}{\phi_j}| \leq
    \normhp{f} \leq \alpha,
  \end{equation*}
  thus $\norminfty{B_{h'}} \leq \alpha$. Since $V_{h'}$ is,
  conditional on $\mf X_n$, a centered Gaussian vector with variance
  $\sigma^2 N_{n,h'}^{-1} \Lba_{h'} \mb X_{h'} \Lba_{h'}$ we have that
  $\mc G_{h'}^{-1} \Lba_{h'} V_{h'}$ is also centered Gaussian, with
  variance
  \begin{equation*}
    \sigma^2 N_{n,h'}^{-1} \mc G_{h'}^{-1} \Lba_{h'} \mb
    X_{h'} \Lba_{h'} \mc G_{h'}^{-1} = \sigma^2 N_{n,h'}^{-1}
    \Lba_{h'}^{-1} \mbt X_{h'}^{-1} \mb X_{h'} \mbt X_{h'}^{-1}
    \Lba_{h'}^{-1}. 
  \end{equation*}
  The variable $\prodsca{\mc G_{h'}^{-1} V_{h'}}{e_1}$ is then
  conditional on $\mf X_n$ centered Gaussian with variance
  \begin{equation*}
    v_{h'}^2 \eqdef \sigma^2 N_{n,h'}^{-1} \prodsca{e_1 \Lba_{h'}^{-1} \mbt
      X_{h'}^{-1} \mb X_{h'} \mbt X_{h'}^{-1} \Lba_{h'}^{-1} }{e_1}
    \leq \sigma^2 N_{n,h'}^{-1} \norm{\Lba_{h'}^{-1}}^2 \norm{\mbt
      X_{h'}^{-1}}^2 \norm{\mb X_{h'}},
  \end{equation*}
  and since clearly $\norm{\mb X_{h'}} \leq \kpa + 1$,
  $\norm{\Lba_{h'}^{-1}} \leq 1$ and $\norm{\mbt X_{h'}^{-1}} =
  \lba^{-1}(\mbt X_{h'}) \leq N_{n,h'}^{1/2}$ we have $v_{h'}^2 \leq
  \sigma^2(\kpa + 1)$ and $\norm{\mc G_{h'}^{-1}} \leq
  \norm{\Lba_{h'}^{-1}} \norm{\mbt X_{h'}^{-1}} \norm{\Lba_{h'}^{-1}}
  \leq N_{n,h'}^{1/2}$.  Finally we have
  \begin{align*}
    |\wh f_{h'}(x_0)| \leq |\prodsca{\mc G_{h'}^{-1} B_{h'}}{e_1}| +
    |\prodsca{\mc G_{h'}^{-1} V_{h'}}{e_1}| &\leq \norm{\mc
      G_{h'}^{-1}} (\norm{B_{h'}} + \sigma \sqrt{\kpa + 1}
    |\gamma_{h'}|) \\
    &\leq \sqrt{\kpa + 1} N_{n,h'}^{1/2} ( \norminfty{B_{h'}} + \sigma
    |\gamma_{h'}| ) \\
    &\leq \sqrt{\kpa + 1} (\alpha \vee 1) N_{n,h}^{1/2} (1 + \sigma
    |\gamma_{h'}| ),
  \end{align*}
  where $\gamma_{h'}$ is, conditional on $\mf X_n$, centered Gaussian
  with variance $v_{h'}^2 \leq 1$. The lemma follows by integrating
  with respect to $\Pfm(\cdot | \mf X_n)$. 
\end{proof}

\begin{proof}[Proof of theorem \ref{thm:NA_adapt_upper_bound}]

  We first work on the event $\{ \wh H_n < H_{n, \omega}^* \}$. By
  definition of $\wh H_n$ we have $\{ \wh H_n < H_{n, \omega}^* \}
  \subset \mc A_{H_{n, \omega}^*}^c$. Uniformly for $f \in \mc
  U(\alpha)$ we have using the lemmas \ref{lem:prob_bandwidth_reject}
  and \ref{lem:bad_case}:
  \begin{align*}
    \Efm \bigl\{ R_{n, \omega}^{-p} | \wh f_n(x_0) &- f(x_0) |^p
    \ind{\wh H_n < H_{n, \omega}^*} | \mf X_n \bigr\} \\
    &\leq (2^p \vee 1) R_{n, \omega}^{-p} \Bigl( \sqrt{ \Efm\{ |\wh
      f_{\wh H_n}(x_0)|^{2p} | \mf X_n \}} + |f(x_0)|^p \Bigr)
    \sqrt{\Pfm\{ \mc A_{H_{n, \omega}^*}^c | \mf X_n \}} \\
    &\leq (2^p \vee 1) \sigma^{-p} (\alpha \vee 1)^p (
    \sqrt{C_{\sigma, 2p, \kpa}} + 1) \sqrt{\kpa + 1} (\log n)^{-p/2} =
    o_n(1). 
  \end{align*}
  Now we work on the event $\{ H_{n, \omega}^* \leq \wh H_n \}$. By
  definition of $\wh H_n$ we have $\{ H_{n, \omega}^* \leq \wh H_n \}
  \subset \mc A_{H_{n, \omega}^*, \wh H_n}$ and using lemma
  \ref{lem:estimator_difference} we get on $\Omega_{H_{n, \omega}^*}$:
  \begin{equation}
    \label{eq:proof_NA_upper_bound_ineq1}
    | \wh f_{\wh H_n}(x_0) - \wh f_{H_{n, \omega}^*}(x_0) | \leq
    C_{p, \kpa, a} \norm{\mc G_{H_{n, \omega}^*}^{-1}} R_{n,
      \omega}. 
  \end{equation}
  Since $s \leq \kpa+1$ we have $k = \ppint{s} \leq \kpa$ and
  $\omega_{f,\kpa}(x_0,h) \leq \omega_{f,k}(x_0,h)$. In view of lemma
  \ref{lem:bias_variance_adapt} and since $f \in \mc F_{H_{n,
      \omega}^*}(x_0,\omega)$ one has on $\Omega_{H_{n, \omega}^*}$:
  \begin{equation*}
    |\wh f_{H_{n, \omega}^*}(x_0) - f(x_0)| \leq \norm{\mc
      G_{H_{n, \omega}^*}^{-1}} \sqrt{\kpa+1}
    \bigl(\omega(H_{n, \omega}^*) + \sigma N_{n,H_{n, \omega}^*}^{-1/2}
    |\gamma_{H_{n, \omega}^*}| \bigr),
  \end{equation*}
  where $\gamma_{H_{n, \omega}^*}$ is, conditional on $\mf X_n$,
  centered Gaussian with $\Efm\{ \gamma_{H_{n, \omega}^*}^2 | \mf X_n
  \} \leq 1$. When $H_{n, \omega}^* < H_{n, \omega}$ we have
  $\omega(H_{n, \omega}^*) \leq \sigma \sqrt{N_{n,H_{n,
        \omega}^*}^{-1} \log n}$. When $H_{n, \omega}^* = H_{n,
    \omega}$ we proceed as in the proof of lemmas
  \ref{lem:prob_bandwidth_reject} and \ref{lem:estimator_difference}
  to prove that
  \begin{equation*}
    \omega(H_{n, \omega}^*) \leq \sigma \sqrt{(1+a) N_{n, H_{n,
          \omega}^*}^{-1} \log n},
  \end{equation*}
  in both cases $\mc H = \mc H_a^{\text{arith}}$ or $\mc H = \mc
  H_a^{\text{geom}}$. Then
  \begin{equation}
    \label{eq:proof_NA_upper_bound_ineq2}
    |\wh f_{H_{n, \omega}^*}(x_0) - f(x_0)| \leq R_{n, \omega} \norm{\mc
      G_{H_{n, \omega}^*}^{-1}} \sqrt{\kpa+1} \bigl(\sqrt{1+a} +
    |\gamma_{H_{n, \omega}^*}| \bigr). 
  \end{equation}
  Finally, the inequalities \eqref{eq:proof_NA_upper_bound_ineq1} and
  \eqref{eq:proof_NA_upper_bound_ineq2} together entail:
  \begin{align*}
    R_{n, \omega}^{-1} |\wh f_n(x_0) - f(x_0)| \ind{ H_{n, \omega}^*
      \leq \wh H_{n, \omega} } \leq \norm{\mc G_{H_{n,
          \omega}^*}^{-1}} (C_{p, \kpa, a} + \sqrt{\kpa+1} (
    \sqrt{1+a} + |\gamma_{H_{n, \omega}^*}|) ),
  \end{align*}
  and the result follows by integration with respect to $\Pfm(\cdot
  |\mf X_n)$. 
\end{proof}

\subsection{Preparatory results and proof of theorem
  \ref{thm:asympt_adaptive_upper_bound}}
\label{sec:preparatory_results_and_proof_of_theorem_2}

Let us denote by $\Pm$ the joint probability of the variables
$(X_i)_{i=1,\ldots,n}$. We define $F_{\nu}(h) \eqdef \int_0^h \nu(t)
dt$. 

\begin{lemma}
  \label{lem:Nh_prob_equiv}
  If $\mu \in \mc R(x_0,\beta)$ one has for any $\von, h > 0$\tup:
  \begin{equation*}
    \forall \von > 0, \quad \Pm \Bigl\{ \Bigl| \frac{N_{n, h}}{2 n
      F_{\nu}(h)} - 1 \Bigr| > \von \Bigr\} \leq 2 \exp \Bigl(
    -\frac{\von^2}{1 + \von / 3} n F_{\nu}(h) \Bigr). 
  \end{equation*}
\end{lemma}
\begin{proof}
  It suffices to use the Bernstein inequality to the sum of
  independent random variables $ Z_i = \ind{|X_i - x_0| \leq h} -
  \Pm\{ |X_1 - x_0| \leq h \}$ for $i=1,\ldots,n$. 
\end{proof}

\begin{lemma}
  \label{lem:rv_r_n_equiv_adapt}
  If $\mu \in \mc R(x_0,\beta)$ for $\beta > -1$, $\omega \in \RV(s)$
  for $s>0$ and $(h_{n, \omega})$ is defined by
  \eqref{eq:bias_variance_det_adapt} then $r_{n, \omega} =
  \omega(h_{n, \omega})$ satisfies
  \begin{equation}
    \label{eq:rv_h_n_asympt_adapt}
    r_{n, \omega} \sim c_{s, \beta} \sigma^{2s / (1 + 2s +
      \beta)} (\log n / n)^{s/(1 + 2s + \beta)} \ell_{\omega, \nu}(
    \log n / n) \text{ as } n \raro + \infty,
  \end{equation}
  where $\ell_{\omega,\nu}$ is slowly varying and $c_{s, \beta} = 2^{s
    / (1 + 2s + \beta)}$. When $\omega(h) = r h^s$ \tup(H{\"o}lder
  regularity\tup) for $r>0$ we have more precisely\tup:
  \begin{equation}
    \label{eq:rv_h_n_asympt_holder_adapt}
    r_{n, \omega} \sim c_{s, \beta} \sigma^{2s / (1 + 2s +
      \beta)} r^{(\beta + 1) / (1 + 2s + \beta)} (\log n / n)^{s / (1
      + 2s + \beta)} \ell_{s,\nu}(\log n / n) \text{ as } n \raro + \infty,
  \end{equation}
  where $\ell_{s,\nu}$ is again slowly varying. 
\end{lemma}

\begin{proof}
  Let us define $G(h) = \omega^2(h) F_{\nu}(h)$. Since $\beta > -1$ we
  have $F_{\nu} \in \RV(\beta+1)$ (see appendix) and $G \in
  \RV(1+2s+\beta)$. The function $G$ is continuous and such that
  $\lim_{h \raro 0^+} G(h) = 0$ in view of
  \eqref{eq:rv_limit_0_or_infty}, since $1+2s+\beta > 0$.  Then for
  $n$ large enough $h_n$ is given by $h_n = G^{\laro}(\sigma^2 \log n
  / 2 n)$ where $G^{\laro}(h) \eqdef \inf\{ y \geq 0 | G(y) \geq h \}$
  is the generalised inverse of $G$. Since $G^{\laro} \in \RV(1/(1 +
  2s + \beta))$ (see appendix) we have $\omega \circ G^{\laro} \in
  \RV(s / (1 + 2s + \beta))$ and we can write $\omega \circ G^{\laro}
  = h^{s / (1 + 2s + \beta)} \ell_{\omega, \nu}(h)$ where
  $\ell_{\omega, \nu}$ is slowly varying.  Thus
  \begin{align*}
    r_n = \omega \circ G^{\laro} \Big( \sigma^2 \frac{\log n}{2 n}
    \Big) &= c_{s, \beta} \sigma^{2s / (1 + 2s + \beta)} (\log n /
    n)^{s/(1 + 2s + \beta)} \ell_{\omega,\nu} \Big( \frac{\sigma^2
      \log n}{2n} \Big) \\
    &\sim c_{s, \beta} \sigma^{2s / (1 + 2s + \beta)} (\log n /
    n)^{s/(1 + 2s + \beta)} \ell_{\omega,\nu} \Big( \frac{\log n}{n}
    \Big) \text{ as } n \raro +\infty,
  \end{align*}
  since $\ell_{\omega, \nu}$ is slowly varying. When $\omega(h) = r
  h^s$ we can write more precisely $h_n = G^{\laro} \big(
  \frac{\sigma^2 \log n}{2 r^2 n} \big)$ where $G(h) = h^{2s}
  F_{\nu}(h)$, so \eqref{eq:rv_h_n_asympt_adapt} and
  \eqref{eq:rv_h_n_asympt_holder_adapt} follow. 
\end{proof}

Let us introduce the following notations: if $\alpha \in \setN$ and
$h>0$ we define
\begin{equation*}
  N_{n,h,\alpha} \eqdef \sum_{|X_i-x_0| \leq h}
  \Bigl(\frac{X_i-x_0}{h}\Bigr)^{\alpha}. 
\end{equation*}
Note that $N_{n,h,0} = N_{n,h}$. For $\von > 0$, we define the event:
\begin{equation*}
  \mrm D_{n,h,\alpha,\von} \eqdef \Bigl\{ \Bigl|
  \frac{N_{n,h,\alpha}}{n F_{\nu}(h)} - C_{\alpha,\beta} \Bigr| \leq
  \von \Bigr\},
\end{equation*}
where $C_{\alpha,\beta}$ is given in section
\ref{sect:regularly_varying_design}. 

\begin{lemma}
  \label{lem:bernstein_adapt}
  For any $\alpha \in \setN$, $\von > 0$ and if $\mu \in \mc R(x_0,
  \beta)$ we have for any positive sequence $(\gamma_n)$ going to $0$
  and when $n$ is large enough\tup:
  \begin{equation}
    \label{eq:bernstein}
    \Pm \bigl\{ \mrm D_{n,\gamma_n,\alpha,\von}^c \bigr\} \leq 2
    \exp \Big( - \frac{\von^2}{ 8 (2 + \von / 3)} n F_{\nu}(\gamma_n)
    \Big). 
  \end{equation}
\end{lemma}

\begin{proof}
  Let us define $Q_{i,n,\alpha} \eqdef \bigl(\frac{X_i-x_0}{\gamma_n}
  \bigr)^{\alpha} \ind{|X_i-x_0| \leq \gamma_n} $, $Z_{i,n,\alpha}
  \eqdef Q_{i,n,\alpha} - \Em\{ Q_{i,n,\alpha} \}$. Since $\mu \in \mc
  R(x_0, \beta)$ one has for $n$ such that
  $[x_0-\gamma_n,x_0+\gamma_n] \subset W$ and $i \in \{ 1,\ldots,n
  \}$:
  \begin{equation*}
    \frac{1}{F_{\nu}(\gamma_n)} \Em \{ Q_{i,n,\alpha} \} =  (1 +
    (-1)^{\alpha}) \frac{\gamma_n^{\beta + 1} \ell_{\nu}(\gamma_n)
    }{\int_0^{\gamma_n} t^{\beta} \ell_{\nu}(t) dt}
    \frac{\int_0^{\gamma_n} t^{\alpha + \beta} \ell_{\nu}(t)
      dt}{\gamma_n^{\alpha + \beta + 1 } \ell_{\nu}(\gamma_n)},
  \end{equation*}
  where $\ell_{\nu}(h) = h^{-\beta}\nu(h)$ is slowly varying (see
  appendix) and in view of \eqref{eq:karamata_equiv} we have
  \begin{equation*}
    \lim_{n \raro +\infty} \frac{1}{F_{\nu}(\gamma_n)} \Em \{
    Q_{i,n,\alpha} \} = C_{\alpha,\beta}. 
  \end{equation*}
  Then for $n$ large enough one has:
  \begin{equation}
    \label{eq:inclusion_lemma_bernstein}
    \Big\{ \bigl| \frac{N_{n,\gamma_n,\alpha}}{n F_{\nu}(\gamma_n)} -
    C_{\alpha,\beta} \bigr| > \von \Big\} \subset \biggl\{ \Bigl|
    \frac{1}{n F_{\nu}(\gamma_n)} \sum_{i=1}^n Z_{i,n,\alpha}\Bigr| >
    \von / 2 \biggr\}. 
  \end{equation}
  We have $ \Em\{{Z_{i,n,\alpha}}\}=0 $, $|Z_{i,n,\alpha}| \leq 2 $. 
  Since
  \begin{equation*}
    b_{n,\alpha}^2 \eqdef \sum_{i=1}^n \Em \{ Z_{i,n,\alpha}^2 \} \leq
    n \Em \{ Q_{i,n,\alpha}^2 \} \leq 2 n F_{\nu}(\gamma_n),
  \end{equation*}
  and the $Z_{i,n,\alpha}$ are independent we can apply Bernstein
  inequality. If $\tau_{n} \eqdef \frac{\von}{2} n F_{\nu}(\gamma_n)$,
  \eqref{eq:inclusion_lemma_bernstein} and Bernstein inequality
  entail:
  \begin{equation*}
    \Pm \bigl\{ \mrm D_{n, \gamma_n, \alpha, \von}^c \bigr\} \leq
    2 \exp \biggl( \frac{-\tau_{n}^2}{ 2(b_{n,\alpha}^2 + 2 \tau_n /3)
    } \biggr) \leq 2 \exp \Big( - \frac{\von^2}{8(2 + \von / 3)} n
    F_{\nu}(\gamma_n) \Big). \qedhere
  \end{equation*}
\end{proof}

Let us introduce for $\von >0$ the event
\begin{equation*}
  \mrm C_{n,\von} \eqdef \{ (1 - \von)h_{n, \omega} < H_{n, \omega}
  \leq (1 + \von) h_{n, \omega} \},
\end{equation*}
where $h_{n, \omega}$ is given by \eqref{eq:bias_variance_det_adapt}. 

\begin{lemma}
  \label{lem:D_int_D_inc_C}
  If $\omega \in \RV(s)$ for $s > 0$ then for any $0 < \von_2 \leq
  1/2$ there exists $0 < \von_3 \leq \von_2$ such that for $n$ large
  enough
  \begin{equation*}
    \mrm D_{n,(1-\von_2)h_{n, \omega}, 0, \von_3} \cap \mrm
    D_{n, (1+ \von_2) h_{n, \omega},
      0, \von_3} \subset \mrm C_{n,\von_2}. 
  \end{equation*}
\end{lemma}

\begin{proof}
  By the definition \eqref{eq:H_n_def_adapt} of $H_{n, \omega}$ we
  have
  \begin{equation*}
    \{ H_{n, \omega} \leq (1 + \von_2) h_{n, \omega} \} = \{
    N_{n, (1 + \von_2) h_{n, \omega}} \geq \sigma^2
    \omega^{-2}((1 + \von_2) h_{n, \omega}) \log n \}. 
  \end{equation*}
  It is clear that $\von_3 \eqdef 1 - (1 - \von_2^2)^{-2} (1 +
  \von_2)^{-2 s} \wedge \von_2 > 0$ for $\von_2$ small enough. We
  recall that $\ell_{\omega}$ stands for the slow term of $\omega$
  (see definition \ref{def:rv_sigma_definition}). Since
  \eqref{eq:slow_variation_def} uniformly over each compact set in
  $(0, +\infty)$ we have when $n$ is large enough that for any $y \in
  [\frac{1}{2}, \frac{3}{2}]$:
  \begin{equation}
    \label{eq:step5_UCT1}
    (1-\von_2^2) \ell_{\omega}(h_{n, \omega}) \leq \ell_{\omega}(
    y h_{n, \omega}) \leq (1 + \von_2^2) \ell_{\omega}(h_{n,
      \omega}),
  \end{equation}
  so \eqref{eq:step5_UCT1} with $y = 1 + \von$ ($\von \leq 1/2$)
  entails in view of \eqref{eq:bias_variance_det_adapt} and since
  $F_{\nu}$ is increasing:
  \begin{align*}
    2(1 - \von_3)n F_{\nu}((1 + \von_2) h_{n, \omega}) &\geq (1 -
    \von_2^2)^{-2} (1 + \von_2)^{-2s} \sigma^2 \omega^{-2}(h_{n,
      \omega}) \log n \\
    &= \sigma^2 \big( (1 + \von_2) h_{n, \omega} \big)^{-2s}
    (1 - \von_2^2)^{-2} \ell_{\omega}^{-2}(h_{n, \omega}) \log n \\
    &\geq \sigma^2 \omega^{-2} ( (1 + \von_2) h_{n, \omega}) \log n. 
  \end{align*}
  Thus
  \begin{equation*}
    \{ N_{n, (1+\von_2) h_{n, \omega}} \geq 2(1 - \von_3) n
    F_{\nu}((1 + \von_2) h_{n, \omega}) \} \subset \{ H_{n,
      \omega} \leq (1 + \von_2) h_{n, \omega} \},
  \end{equation*}
  and similarly on the other side we have for $n$ large enough
  \begin{equation*}
    \{ N_{n, (1 - \von_2) h_{n, \omega}} \leq 2(1 + \von_3)n
    F_{\nu}((1 - \von_2) h_{n, \omega}) \}
    \subset \{ (1 - \von_2) h_{n, \omega} < H_{n, \omega} \},
  \end{equation*}
  thus the lemma. 
\end{proof}

Let us denote $\mc G_n \eqdef \mc G_{H_{n, \omega}}$ and introduce the
events $\mrm A_{n,\von} \eqdef \bigl\{ | \lba(\mc G_n) - \lba_{\kpa,
  \beta} | \leq \von \bigr\}$ for $\von > 0$ and for $\alpha \in
\setN$
\begin{equation*}
  \mrm B_{n, \alpha, \von} \eqdef \Bigl\{ \Bigl| \frac{1}{n F_{\nu}(h_n)}
  \sum_{|X_i - x_0| \leq H_n} \Bigl( \frac{X_i - x_0}{h_n}
  \Bigr)^{\alpha} - C_{\alpha,\beta} \Bigr| \leq \von \Bigr\}. 
\end{equation*}

\begin{lemma}
  \label{lem:A_n_von_subset}
  If $\omega \in \RV(s)$ for $s > 0$ and $\mu \in \mc R(x_0, \beta)$
  for $\beta > -1$ we can find for any $0 < \von \leq \frac{1}{2}$ an
  event $\mc A_{n,\von} \in \mf X_n$ such that for $n$ large enough
  \begin{equation}
    \label{eq:A_n_von_subset_things}
    \mc A_{n,\von} \subset \mrm A_{n, \von} \cap \mrm B_{n, 0, \von}
    \cap \mrm C_{n, \von},
  \end{equation}
  and
  \begin{equation}
    \label{eq:A_n_von_deviation}
    \Pm \{ \mc A_{n, \von}^c \} \leq 4(\kpa + 2) \exp \big( -c_{\beta,
      \sigma, \von} r_n^{-2 }\big). 
  \end{equation}
\end{lemma}

\begin{proof}
  Using the fact that $\lba(M) = \inf_{\norm{x} = 1} \prodsca{x}{Mx}$
  for any symmetrical matrix $M$ and since $\mc G_n$ and $\mc G$ are
  symmetrical we get
  \begin{equation*}
    \bigcap_{\alpha=0}^{2 \kpa} \Bigl\{ \bigl| (\mc G_n)_{j,l} - (\mc
    G)_{j,l} \bigr| \leq \frac{\von}{(1 + \kpa)^2} \Bigr\} \subset
    \mrm A_{n,\von}. 
  \end{equation*}
  Since $|(\mc G)_{j,l}| \leq 1$ we can find easily $0 < \von_1 \leq
  \von$ such that for any $0 \leq j,l \leq \kpa$
  \begin{equation*}
    \mrm B_{n, j + l, \von_1} \cap \mrm B_{n, 2j, \von_1} \cap \mrm
    B_{n, 2l, \von_1} \subset \Bigl\{ \bigl| (\mc G_n)_{j,l} - (\mc
    G)_{j,l} \bigr| \leq \frac{\von}{(1 + \kpa)^2} \Bigr\},
  \end{equation*}
  and then
  \begin{equation*}
    \bigcap_{\alpha=0}^{2 \kpa} \mrm B_{n, \alpha, \von_1} \subset \mrm
    A_{n, \von}. 
  \end{equation*}
  We define $ \von_2 \eqdef \frac{2^{\kpa}}{5 \times 3^{\kpa}} \von_1$
  and $\von_3$ such that $\frac{(2 + \von_3)(1 + \von_3)^{\beta +
      2}}{2 - \von_3} = 1 + \von_2$. Since $h \mapsto N_{n,h}$ is
  increasing we have
  \begin{equation*}
    \mrm C_{n, \von_3} \subset \{ N_{n,(1 - \von_3)h_n} \leq N_{n,H_n}
    \leq N_{n,(1 + \von_3) h_n} \},
  \end{equation*}
  and using lemma \ref{lem:D_int_D_inc_C} we can find $\von_4 \leq
  \von_3$ such that
  \begin{equation*}
    \mrm D_{n, (1 - \von_3) h_n, 0, \von_4} \cap \mrm D_{n, (1 +
      \von_3) h_n, 0, \von_4} \subset \mrm C_{n,\von_3}. 
  \end{equation*}
  In view of \eqref{eq:slow_variation_def} and since $\ell_{\nu}(h)
  \eqdef F_{\nu}(h) h^{-(\beta + 1)}$ is slowly varying we have for
  $n$ large enough and any $0 < \von_3 \leq 1/2$
  \begin{equation}
    \label{eq:lets_use_UCT_one_more_time}
    \ell_{\nu}((1 + \von_3) h_n) \leq (1
    + \von_3) \ell_{\nu}(h_n) \text{ and } \ell_{\nu}((1 - \von_3) h_n) \geq (1
    - \von_3)\ell_{\nu}(h_n),
  \end{equation}
  thus
  \begin{equation*}
    \mrm D_{n, (1 - \von_3) h_n, 0, \von_4} \cap
    \mrm D_{n, (1+\von_3) h_n, 0, \von_4} \cap \mrm D_{n, h_n, 0, \von_3}
    \subset \mrm E_{n, \von_2} \eqdef \Bigl\{ \Bigl| \frac{N_{n, H_n}}{N_{n,
        h_n}} - 1 \Bigr| \leq \von_2 \Bigr\},
  \end{equation*}
  and on $\mrm D_{n, (1 - \von_3)h_n, 0, \von_4} \cap \mrm D_{n, (1 +
    \von_3)h_n, 0, \von_4} \cap \mrm D_{n, h_n, 0, \von_3}$ we have
  \begin{align*}
    \frac{1}{n F_{\nu}(h_n)} \Bigl| \sum_{|X_i - x_0| \leq H_n} \Bigl(
    \frac{X_i - x_0}{h_n} \Bigr)^{\alpha} - N_{n,h_n,\alpha} \Bigr|
    &\leq \Bigl( \frac{H_n \vee h_n}{h_n} \Bigr)^{\alpha}
    \frac{N_{n,h_n}}{n F_{\nu}(h_n)} \Bigl|
    \frac{N_{n,H_n}}{N_{n,h_n}} - 1 \Bigr| \\
    &\leq (1 + \von_3)^{\alpha} (2 + \von_3) \von_2 \leq \von_1 / 2,
  \end{align*}
  since $\von_3 \leq 1/2$. Then we have since $\von_4 \leq \von_3 \leq
  \von_2 \leq \frac{\von_1}{2}$
  \begin{align*}
    \mrm D_{n,(1-\von_3)h_n, 0, \von_4} \cap \mrm D_{n,(1+\von_3) h_n,
      0, \von_4} \cap \mrm D_{n,h_n,0,\von_4} \cap \mrm D_{n, h_n,
      \alpha, \von_4} \subset \mrm B_{n, \alpha, \von_1},
  \end{align*}
  and finally
  \begin{equation*}
    \mc A_{n, \von} \eqdef \mrm D_{n, (1 - \von_3) h_n, 0, \von_4}
    \cap \mrm D_{n, (1 + \von_3) h_n, 0, \von_4} \cap \mrm D_{n, h_n,
      0, \von_4} \cap \bigcap_{\alpha = 0}^{2\kpa} \mrm D_{n, h_n,
      \alpha, \von_4} \subset \mrm A_{n,\von} \cap \mrm B_{n, 0, \von}
    \cap \mrm C_{n, \von},
  \end{equation*}
  thus \eqref{eq:A_n_von_subset_things}. Using lemma
  \ref{lem:bernstein_adapt} we obtain easily in view of
  \eqref{eq:lets_use_UCT_one_more_time} and
  \eqref{eq:bias_variance_det_adapt} for $n$ large enough
  \begin{equation*}
    \Pm \{ \mc A_{n, \von}^c \} \leq 4(\kpa + 2) \exp \Big( \frac{
      -\von_4^2}{ 4 (2 + \von_4 / 3) } 2^{-(\beta + 2)} \sigma^2
    r_n^{-2} \log n \Big),
  \end{equation*}
  thus \eqref{eq:A_n_von_deviation} and the lemma follows. 
\end{proof}

\begin{proof}[Proof of theorem \ref{thm:asympt_adaptive_upper_bound}]
  
  Since $\mc H = \mc H_{1}^{\tup{arith}}$ we have $H_{n, \omega} =
  H_{n, \omega}^*$ and $\lba_{n, \omega} = \lba(\mc G_{H_{n,
      \omega}})$. We can assume without generality loss that $\von
  \eqdef \varrho - 1 \leq \frac{1}{2} \wedge \lba_{\kpa, \beta}$. We
  consider the event $\mc A_{n, \von}$ from lemma
  \ref{lem:A_n_von_subset}.  Clearly, we have for $n$ large enough
  $\mc A_{n,\von} \subset \Omega_{H_{n, \omega}}$ and $\mc F_{\varrho
    h_{n, \omega}}(x_0, \omega) \subset \mc F_{H_{n, \omega}}(x_0,
  \omega)$. In view of \eqref{eq:A_n_von_subset_things} and theorem
  \ref{thm:NA_adapt_upper_bound} we have uniformly for $f \in \Sigma$:
  \begin{align*}
    \Efm \{ r_n^{-p} |\wh f_n(x_0) - f(x_0)|^p \ind{\mc A_{n,\von}} \}
    &\leq (1 - \von)^{-p/2} \Efm\{ R_n^{-p} |\wh f_n(x_0) -
    f(x_0)|^p \ind{\Omega_{H_n}} \} \\
    &\leq (1 - \von)^{-p/2} c_1 (\lba_{\kpa,\beta} - \von)^{-p} (1 +
    o_n(1)). 
  \end{align*}
  Now we work on the complementary $\mc A_{n,\von}^c$. Using lemma
  \ref{lem:bad_case} and equation \eqref{eq:A_n_von_deviation} we get
  since $f \in \mc U(\alpha)$ and $N_{n, h} \leq n$:
  \begin{align*}
    \Efm \{ r_n^{-p} |\wh f_n(x_0) - f(x_0)|^p \ind{\mc A_{n,\von}^c}
    \} &\leq (2^p \vee 1) r_n^{-p} \big( \sqrt{ \Efm\{ |\wh
      f_n(x_0)|^{2 p} \}} + \alpha^p \big) \sqrt{ \Pm\{ \mc A_{n,
        \von}^c \}} \\
    &\leq (2^p \vee 1) (\alpha \vee 1)^p (\sqrt{C_{\sigma, 2p, \kpa}}
    + 1) n^{p/2} r_n^{-p} \sqrt{ \Pm\{ \mc A_{n, \von}^c \}} = o_n(1),
  \end{align*}  
  thus we have proved \eqref{eq:adapt_upper_bound_asympt} and
  \eqref{eq:upper_bound_r_n_equiv_adapt} follows from lemma
  \ref{lem:rv_r_n_equiv_adapt}. 
\end{proof}

\subsection{Computation of the example}
\label{sec:computation_exemple}

\begin{lemma}
  \label{lem:lambert_inversion}
  Let $a \in \setR$ and $b > 0$. If $G(h) = h^{b} (\log(1/h))^{a}$,
  then we have
  \begin{equation*}
    G^{\laro} (h) \sim b^{a / b} h^{1/ b} (\log(1 / h))^{- a / b }
    \text{ as } h \raro 0^+. 
  \end{equation*}
\end{lemma}
The proof of this lemma can be found in Ga{\"\i}ffas
(2004)\nocite{gaiffas04a}. Using this lemma, we obtain that an
equivalent of $h_n$ (see \eqref{eq:bias_variance_det_adapt}) is
\begin{equation*}
  (1 + 2s + \beta)^{(\alpha + 2 \gamma) / (1 + 2s + \beta)} \Bigl(
  \frac{\sigma}{r} \Big)^{2 / (1 + 2s + \beta)} \big( 2 n (\log n)^{
    \alpha + 2 \gamma - 1 } \big)^{-1 / (1 + 2s + \beta)},
\end{equation*}
and since $\omega(h) = r h^s (\log(1 / h))^{\gamma}$ we find that an
equivalent of $r_n$ (up to a constant depending on $s, \beta, \gamma,
\alpha$) is \eqref{eq:convergence_rate_exemple}.

\subsection{Proof of the lower bound}
\label{sec:proof_lower_bound}

The proof of theorem \ref{thm:lower_bound} is similar to the proof of
theorem 3 in Brown and Low (1996)\nocite{brown_low_constrained96}. It
is based on the next theorem which can be found in Cai \etal
(2004)\nocite{cai_low_zhao}. This result is a general constrained risk
inequality and is very useful for several statistical problems, for
example superefficiency, adaptation and so on.

Let $p > 1$ and $q$ be such that $\frac{1}{p} + \frac{1}{q} = 1$ and
$X$ be a real random variable having distribution $\Prob_{\tta}$ with
density $f_{\tta}$ with respect to some measure $m$. The parameter
$\tta$ can take two values $\tta_1$ or $\tta_2$. We want to estimate
$\tta$ based on $X$. For any estimator $\delta$ based on $X$ we define
its risk by
\begin{equation*}
  R_p(\delta, \tta) \eqdef \mbb E_{\tta} \{ |\delta(X) - \tta |^p \}. 
\end{equation*}
We define $s(x) = f_{\tta_2}(x) / f_{\tta_1}(x)$ and $\Delta = |\tta_2
- \tta_1|$. Let
\begin{equation*}
  I_q = I_q(\tta_1, \tta_2) \eqdef \big( \mbb E_{\tta_1} \{ s^q(X) \}
  \big)^{1 / q}. 
\end{equation*}

\begin{theorem}[Cai, Low and Zhao (2004)]
  \label{thm:constrained_risk_inequality}
  If $\delta$ is such that $R_p(\delta, \tta_1) \leq \von^p$ and if
  $\Delta > \von I_q$ we have
  \begin{equation*}
    R_p(\delta, \tta_2) \geq (\Delta - \von I_q )^p \geq \Delta^p
    \Big(1 - \frac{p \von I_q}{\Delta} \Big). 
  \end{equation*}
\end{theorem}

\begin{proof}[Proof of theorem \ref{thm:lower_bound}]
  Since $\limsup_n \psi_n^{-p} n^{\gamma p} \mbb E_{f_0, \mu} \{ | \wh
  f_n(x_0) - f_0(x_0) |^p \} = C < \infty$ we have for $n \geq N$
  \begin{equation*}
    \mbb E_{f_0, \mu} \{ | \wh f_n(x_0) - f_0(x_0) |^p \} \leq 2 C
    \psi_n^p n^{-\gamma p}. 
  \end{equation*}
  Let $g$ be $k$ times differentiable with support included in $[-1,
  1]$, $g(0) > 0$ and such that for any $|x| \leq \delta$, $|
  g^{(k)}(x) - g^{(k)}(0) | \leq k!  |x|^{s-k}$. Such a function
  clearly exists. We define
  \begin{equation*}
    f_1(x) \eqdef f_0(x) + (r - r') \rho_n^s g \Big( \frac{x -
      x_0}{\rho_n} \Big),
  \end{equation*}
  where $\rho_n$ is the smallest solution to
  \begin{equation*}
    r h^{s} = \sigma \sqrt{ \frac{b \log n}{2 n F_{\nu}(h)}},
  \end{equation*}
  where $b = 2 g_{\infty}^{-2} (p - 1) \gamma$ and $g_{\infty} \eqdef
  \sup_x |g(x)|$. We clearly have $f_1 \in \mc F_{\delta}(x_0,
  \omega)$. Let $\mbb P_0^n, \mbb P_1^n$ be the joint laws of the
  observations \eqref{eq:regression_model} when respectively $f = f_0$
  or $f = f_1$. A sufficient statistic for $\{ \mbb P_0^n, \mbb P_1^n
  \}$ is given by $\displaystyle T_n \eqdef \log \frac{\mrm d \mbb
    P_0^n}{\mrm d \mbb P_1^n}$, and
  \begin{equation*}
    T_n \sim
    \begin{cases}
      \displaystyle
      \mc N \big( - \frac{v_n}{2}, v_n \big) &\text{ under } \mbb P_0^n, \\
      \displaystyle \mc N \big( \frac{v_n}{2}, v_n \big) &\text{ under
      } \mbb P_1^n,
    \end{cases}
  \end{equation*}
  where 
  \begin{equation*}
    v_n = \frac{n}{\sigma^2} \norm{f_0 - f_1}_{L^2(\mu)}^2 =
    \frac{n}{\sigma^2} \int ( f_0(x) - f_1(x) )^2 \mu(x) dx \leq 2 (p
    - 1) \gamma \log n. 
  \end{equation*}
  An easy computation gives $I_q = \exp( \frac{ v_n (q - 1) }{2} )
  \leq n^{\gamma}$ thus taking $\delta_n = \wh f_n(x_0)$, $\tta_2 =
  f_1(x_0)$, $\tta_1 = f_0(x_0)$ and $\von = \psi_n$ entails using
  theorem \ref{thm:constrained_risk_inequality}
  \begin{equation*}
    R_p(\delta_n, \tta_2) \geq \big( (r - r') \rho_n^s g(0) - 2 C \psi_n
    n^{ -\gamma} n^{\gamma} \big)^p \geq (r - r')^p \rho_n^{s p} g^p(0) (1
    - o_n(1)),
  \end{equation*}
  since $\lim_n \psi_n / \rho_n^s \raro 0$, and the theorem follows. 
\end{proof}

\appendix

\section{Some facts on regular variation}
\label{sect:regular_variation_tech}

We recall here briefly some results about regularly varying functions.
The results stated in this section can be found in Senata
(1976)\nocite{senata76}, Geluk and de Haan
(1987)\nocite{geluk_de_haan87} and Bingham \etal (1989)\nocite{bgt89}.

Let $\ell$ be in all the following a slowly varying function. An
important result is that the property
\begin{equation}
  \label{eq:slow_variation_def}
  \lim_{h \raro 0^+} \ell(y h) / \ell(h) = 1
\end{equation}
actually holds \emph{uniformly} for $y$ in any compact set of $(0,
+\infty)$. If $R \in \RV(\alpha_1)$ and $R \in \RV(\alpha_2)$ we have
\begin{itemize}
\item $R_1 \times R_2 \in \RV(\alpha_1 + \alpha_2)$,
\item $R_1 \circ R_2 \in \RV(\alpha_1 \times \alpha_2)$. 
\end{itemize}
If $R \in \RV(\gamma)$ with $\gamma \in \setR - \{ 0 \}$ then as $h
\raro 0^+$ we have
\begin{equation}
  \label{eq:rv_limit_0_or_infty}
  R(h) \raro
  \begin{cases}
    0& \text{ if } \gamma > 0, \\
    +\infty& \text{ if } \gamma <0. 
  \end{cases}
\end{equation}
If $\gamma > -1$, one has:
\begin{equation}
  \label{eq:karamata_equiv}
  \int_{0}^{h} t^{\gamma} \ell(t) dt \sim (1 + \gamma)^{-1}
  h^{1+\gamma} \ell(h) \text{ as } h \raro 0^+,
\end{equation}
and then $h \mapsto \int_0^h t^{\gamma} \ell(t) dt$ is regularly
varying of index $1+\gamma$. This result is known as the Karamata
theorem. If $R$ is continuous we define the generalised inverse as
\begin{equation*}
  R^{\laro}(y) = \inf\{ h \geq 0 \text{ such that } R(h) \geq y \}. 
\end{equation*}
If $R \in \RV(\gamma)$ for some $\gamma > 0$ then there exists $R^{-}
\in \RV(1/\gamma)$ such that
\begin{equation}
  \label{eq:asympt_inverve}
  R ( R^{-} (h) ) \sim R^{-} ( R(h) ) \sim h \text{ as } h \raro
  0^+,
\end{equation}
and $R^{-}$ is unique up to an asymptotic equivalence. Moreover, one
version of $R^{-}$ is $R^{\laro}$.


\subsection*{Acknowledgements}

I wish to thank my adviser Marc Hoffmann for helpful advices and
encouragements.

\bibliographystyle{acm} 
\bibliography{biblio}

\end{document}